\documentclass[a4paper,onecolumn, superscriptaddress,10pt, accepted=2022-07-11, issue=3, volume=3, shorttitle=papers/compositionality-3-3]{compositionalityarticle}
\pdfoutput=1
\usepackage[utf8]{inputenc}
\usepackage[english]{babel}
\usepackage[T1]{fontenc}
\usepackage{float}
\usepackage{amssymb,amsmath,stmaryrd,microtype,xr,enumitem,accents}
\usepackage[numbers,sort&compress]{natbib}
\usepackage[all,cmtip]{xy}

\newenvironment{proof}{\noindent\textbf{Proof: }}{\hfill$\square$\\}

\keywords{enriched semicategory,nonunital monoidal structure,combinatorial model category, Quillen equivalence, locally presentable category, topologically enriched category}

\newcommand{\C}{\mathcal{C}}
\newcommand{\D}{\mathcal{D}}
\newcommand{\K}{\mathcal{K}}

\newcommand{\de}{\partial}
\newcommand{\p}{\times}
\renewcommand{\vec}{\overrightarrow}
\renewcommand{\P}{\mathbb{P}}

\newtheorem{thm}{Theorem}[section]
\newtheorem{prop}[thm]{Proposition}
\newtheorem{lem}[thm]{Lemma}
\newtheorem{exa}[thm]{Example}
\newtheorem{cor}[thm]{Corollary}

\newtheorem{defn}[thm]{Definition}
\newtheorem{nota}[thm]{Notation}

\newtheorem{defnot}[thm]{Definition and notation}
\newcommand{\bd}{\begin{defn}}
\newcommand{\ed}{\end{defn}}
\newcommand{\bdn}{\begin{defnot}}
\newcommand{\edn}{\end{defnot}}
\newcommand{\bp}{\begin{prop}}
\newcommand{\ep}{\end{prop}}
\newcommand{\bth}{\begin{thm}}
\renewcommand{\eth}{\end{thm}}
\newcommand{\bpf}{\begin{proof}}
\newcommand{\epf}{\end{proof}}
\newcommand{\bc}{\begin{cor}}
\newcommand{\ec}{\end{cor}}

\newcommand{\fL}[1]{\ar@{->}[ll]_-{#1}}
\newcommand{\fR}[1]{\ar@{->}[rr]^-{#1}}
\newcommand{\fRr}[1]{\ar@{->}[rrr]^-{#1}}
\newcommand{\fD}[1]{\ar@{->}[dd]_-{#1}}
\newcommand{\fU}[1]{\ar@{->}[uu]^-{#1}}
\newcommand{\f}[2]{\ar@{->}[#1]|{#2}}
\newcommand{\ff}[2]{\ar@2{->}[#1]|{#2}}
\newcommand{\frr}[1]{\ar@{->}[rrrr]^-{#1}}

\newcommand{\fl}[1]{\ar@{->}[l]_-{#1}}
\newcommand{\fr}[1]{\ar@{->}[r]^-{#1}}
\newcommand{\fd}[1]{\ar@{->}[d]_-{#1}}
\newcommand{\fu}[1]{\ar@{->}[u]^-{#1}}
\newcommand{\arrow}[2]{{\begin{array}{c}
		#1\\
		\downarrow\\
		#2
\end{array}}}

\renewcommand{\top}{{\mathbf{Top}}}

\newcommand{\iso}{\cong}

\newcommand{\vI}{\vec{I}}
\renewcommand{\leq}{\leqslant}
\renewcommand{\geq}{\geqslant}

\newcommand{\moore}{{\mathbb{M}}}
\newcommand{\lmoore}{\mathbb{M}_!}

\newcommand{\topdgr}{[\mathcal{P}^{op},\top]}

\def\cartesien{%
  \ar@{-}[]+R+<6pt,-2pt>;[]+RD+<6pt,-6pt>%
  \ar@{-}[]+D+<2pt,-6pt>;[]+RD+<6pt,-6pt>%
}
\def\cocartesien{%
  \ar@{-}[]+L+<-6pt,+2pt>;[]+LU+<-6pt,+6pt>%
  \ar@{-}[]+U+<-2pt,+6pt>;[]+LU+<-6pt,+6pt>%
}
\def\hocartesien{%
  \ar@{-}[]+R+<6pt,-2pt>;[]+RD+<6pt,-6pt>_{h}%
  \ar@{-}[]+D+<2pt,-6pt>;[]+RD+<6pt,-6pt>%
}
\def\hococartesien{%
  \ar@{-}[]+L+<-6pt,+2pt>;[]+LU+<-6pt,+6pt>_{h}%
  \ar@{-}[]+U+<-2pt,+6pt>;[]+LU+<-6pt,+6pt>%
}

\newcommand{\brm}[1]{\rm{\mathbf{#1}}}

\newcommand{\dtop}{{\brm{Flow}}}
\newcommand{\dtopP}{{\mathcal{P}\brm{Flow}}}
\newcommand{\dtopG}{{\mathcal{G}\brm{Flow}}}

\newcommand{\set}{{\brm{Set}}}

\newcommand{\ttop}{{\brm{TOP}}}

\newcommand{\glob}{{\rm{Glob}}}

\newcommand{\globP}{{\rm{Glob}}}

\DeclareMathOperator{\id}{Id}

\DeclareMathOperator{\Obj}{Obj}
\DeclareMathOperator{\Mor}{Mor}

\newcommand{\liminj}{\varinjlim}
\newcommand{\limproj}{\varprojlim}

\newcommand{\Dcat}{{\mathbf{CAT}}}

\newcommand{\ptop}[1]{{\brm{{#1}dTop}}}

\makeatletter
\def\varholim@#1#2{%
  \vtop{\m@th\ialign{##\cr
    \hfil$#1\operator@font holim$\hfil\cr
    \noalign{\nointerlineskip\kern1.5\ex@}#2\cr
    \noalign{\nointerlineskip\kern-\ex@}\cr}}%
}
\def\holimproj{%
  \mathop{\mathpalette\varholim@{\leftarrowfill@\textstyle}}\nmlimits@
}
\def\holiminj{%
  \mathop{\mathpalette\varholim@{\rightarrowfill@\textstyle}}\nmlimits@
}
\makeatother

\makeatletter
\renewcommand*{\@opargbegintheorem}[3]{\trivlist
	\item[\hskip \labelsep{\bfseries #1\ #2}] \textbf{(#3)}\ \itshape}
\makeatother

\DeclareMathOperator{\cell}{{\brm{cell}}}
\DeclareMathOperator{\cof}{{\brm{cof}}}
\DeclareMathOperator{\inj}{{\brm{inj}}}

\DeclareMathOperator{\cocyl}{{Path}}
\newcommand{\adj}[4]{\xymatrix@1{{#1}\ar@/^0.8em/[r]^-{#2} \ar@{}[r]|-{\perp} & \ar@/^0.8em/[l]^-{#3} {#4}}}
\DeclareMathOperator{\CC}{CC}
\newcommand{\thin}[1]{{\pi_0}(#1)}
\newcommand{\ot}{\otimes}

\begin{document}

\title{Homotopy theory of Moore flows (I)}
\date{}
\author{Philippe Gaucher}
\orcid{0000-0003-0287-6252}
\homepage{http://www.irif.fr/{\~{}}gaucher}
\affiliation{Universit\'e Paris Cit\'e, CNRS, IRIF, F-75013, Paris, France}

\maketitle

\begin{abstract}

\emph{Erratum, 11 July 2022: This  is an updated version of the original paper \cite{Moore1} in which the notion of reparametrization category was incorrectly axiomatized. Details on the changes to the original paper are provided in the Appendix.}

	A reparametrization category is a small topologically enriched semimonoidal category such that the semimonoidal structure induces a structure of a semigroup on objects, such that all spaces of maps are contractible and such that each map can be decomposed (not necessarily in a unique way) as a tensor product of two maps. A Moore flow is a small semicategory enriched over the biclosed semimonoidal category of enriched presheaves over a reparametrization category. We construct the q-model category of Moore flows. It is proved that it is Quillen equivalent to the q-model category of flows. This result is the first step to establish a zig-zag of Quillen equivalences between the q-model structure of multipointed $d$-spaces and the q-model structure of flows. 
\end{abstract}

\tableofcontents

\section{Introduction}

\subsection*{Presentation}

The q-model category~\footnote{We use the terminology of \cite{ParamHomTtheory}. The h-model structure and the m-model structure are introduced in \cite{QHMmodel}.} of flows (i.e. small semicategories or small nonunital categories enriched over topological spaces) $\dtop$ was introduced in \cite{model3}. The motivation is the study of concurrent processes up to homotopy. The identities are removed from the structure to ensure the functoriality of some constructions. There are many other topological models of concurrency \cite{mg,MR2545830,equilogical,distinguishedcube1,distinguishedcube2} (the list does not pretend to be exhaustive). They are all of them based on the same kind of ideas. They contain a set of states with a notion of execution path and of homotopy between them to model concurrency. They are introduced for various purposes, some of them having a specific mathematical behaviour. The notion of Grandis' $d$-space \cite{mg} seems to be very successful, probably because of its simplicity. A $d$-space is just a topological space equipped with a distinguished set of continuous maps playing the role of the execution paths of a concurrent process. The distinguished set of continuous maps playing the role of execution paths satisfies various natural axioms. On the contrary, flows are difficult to use because they are roughly speaking $d$-spaces without an underlying topological space. They just have an underlying homotopy type as defined in \cite{4eme}. This fact makes the category of flows quite unwieldy, although it has mathematical features not shared by the other topological models of concurrency, like the possibility of defining functorial homology theories for branching and merging areas of execution paths \cite{exbranch,3eme}. 

An attempt to build a bridge between the $d$-spaces of \cite{mg} and the flows of \cite{model3} was made in \cite{mdtop} by introducing the q-model category $\ptop{\mathcal{G}}$ of multipointed $d$-spaces. The latter are a variant of Grandis' notion of $d$-space introduced in \cite{mg} which is as closed as possible to the notion of flow. The most important feature of this notion is that multipointed $d$-spaces do have an underlying topological space.

However, we had to face an unexpected problem which could be summarized as follows. There exists a ``categorization'' functor $cat:\ptop{\mathcal{G}}\to\dtop$ such that the total left derived functor in the sense of \cite{HomotopicalCategory} induces an equivalence of categories between the homotopy categories of the q-model structure of $\ptop{\mathcal{G}}$ and the q-model structure of $\dtop$ \cite[Theorem~7.5]{mdtop}. This functor was a plausible candidate for a Quillen equivalence. Unfortunately, it is neither a left adjoint nor a right adjoint. The origin of the problem is that the composition of paths in multipointed $d$-spaces is associative up to homotopy whereas the composition of paths in flows is strictly associative. 

The solution proposed in this paper to overcome the above problem is to introduce the category $\dtopG$ of \textit{$\mathcal{G}$-flows}, also called \textit{Moore flows} in the title of the paper (cf. Definition~\ref{MooreFlow} and Definition~\ref{paramG}). Roughly speaking, $\mathcal{G}$-flows are flows with morphisms having a length which is a strictly positive real number. Such an object also contains the information about the way for execution paths to be reparametrized. The main result of this paper can be summarized in the following theorem:

\bth (Theorem~\ref{eq1})
The q-model category of $\mathcal{P}$-flows for any reparametrization category $\mathcal{P}$, and in particular the q-model category of Moore flows, and the q-model category of flows are Quillen equivalent.
\eth

A companion paper \cite{Moore2} proves that the q-model category of $\mathcal{G}$-flows and the q-model category of multipointed $d$-spaces are  Quillen equivalent. It also supplies a new and more conceptual proof of \cite[Theorem~7.5]{mdtop}. Note that unlike what I wrote in the motivation section of \cite{dgrtop}, the results of this paper and of the companion paper \cite{Moore2} are independent of \cite{model2}.

The notion of \textit{reparametrization category} of Definition~\ref{def-reparam} seems to be a new method to deal with Moore paths. The specific interest of this notion in directed homotopy theory is to have the terminal category and the two categories $\mathcal{G}$ of Proposition~\ref{paramG} and $\mathcal{M}$ of Proposition~\ref{paramM} in the same formalism. Indeed, Grandis' $d$-spaces, besides being not multipointed and having identity execution paths, have a set of execution paths invariant by the action of the monoid $\mathcal{M}(1,1)$ of nondecreasing surjective maps from $[0,1]$ to itself in their definition. On the contrary, the notion of multipointed $d$-space used in the companion paper \cite{Moore2} has a set of execution paths invariant by the action of the group $\mathcal{G}(1,1)$ of nondecreasing homeomorphisms from $[0,1]$ to itself only, which is sufficient for computer scientific formalization. The reason of the latter choice, which dates back to \cite{mdtop}, is that it makes the cellular objects of the q-model structure of multipointed $d$-spaces much simpler to understand because their execution paths have no stop point. The presence of stop points in execution paths creates technical complications indeed. This problem is studied in \cite{MR2369163} for the continuous paths of a Hausdorff topological space (the continuous paths without stop point are called \textit{regular}). The formalism of reparametrization category will enable us to treat the case of $\mathcal{M}$ in subsequent papers. Note that unlike the classifying space $B\mathcal{G}(1,1)$, the classifying space $B\mathcal{M}(1,1)$ is contractible by a result of Lawson \cite[Theorem~8.1]{dgrtop}. Therefore, in theory, it should be possible to drop somehow the enriched structure in the definition of a $\mathcal{M}$-flow as given in Definition~\ref{PF}. 

The notion of $\mathcal{P}$-flow of Definition~\ref{PF} for a general reparametrization category $\mathcal{P}$ is primarily designed to better understand the topology of the spaces of execution paths in a concurrent process. This will become more apparent in the companion paper \cite{Moore2} in which a lot of new results about the topology of the spaces of execution paths are expounded in the globular setting, i.e. for cellular objects of the q-model structure of multipointed $d$-spaces. For the precubical setting, the reader will refer to the works of Paliga, Raussen and Ziemia\'nski \cite{MR2521708,MR3722069,MR4070250,raussen2021strictifying,paliga2021configuration}. 

The right Quillen functor $\moore^{\mathcal{G}}$ from multipointed $d$-spaces to $\mathcal{G}$-flows defined in the companion paper consists of forgetting the underlying topological space. From a topological space $X$, on can obtain a multipointed $d$-space $\widetilde{X}=(X,X^0,\P X)$ with $X^0=X$ and where $\P X$ is the set of \textit{all} continuous maps from $[0,1]$ to $X$. The $\mathcal{G}$-flow $\moore^{\mathcal{G}}(\widetilde{X})$ is then an object which contains the same information as the so-called Moore path (semi)category of the topological space $X$. It means that the notion of $\mathcal{P}$-flow is an abstraction of the Moore path (semi)category of a topological space. And for this reason, one of the anonymous referees noticed that there could be a connection with the papers \cite{MR3903058,north2019typetheoretic} to be investigated, the first one introducing a new family of models of type theory based on Moore paths, the second one abstracting the notion of Moore path to characterize type-theoretic weak factorization systems.

Throughout the paper, a general reparametrization category $\mathcal{P}$ in the sense of Definition~\ref{def-reparam} is used. The use of the category of  reparametrization $\mathcal{G}$ of Proposition~\ref{paramG} will be required only in the companion paper \cite{Moore2}.

\subsection*{Outline of the paper}

\begin{itemize}[leftmargin=*]
\item Section~\ref{reminder-space} is a short reminder about $\Delta$-generated spaces and their q-model structure. 
\item Section~\ref{presheaf-over-enriched} is a short reminder about enriched presheaves of topological spaces. Some notations are stated and some formulae are recalled and named for helping to explain the subsequent calculations.
\item Section~\ref{param} introduces the notion of \textit{reparametrization category}. We give two examples of such a category (Proposition~\ref{paramG} and Proposition~\ref{paramM}). The notion of reparametrization category given here is very general. 
\item Section~\ref{presheaf-over-param} sketches the theory of enriched presheaves over a reparametrization category $\mathcal{P}$, called \textit{$\mathcal{P}$-spaces}. The main result of Section~\ref{presheaf-over-param} is that it is possible to define a structure of \textit{biclosed semimonoidal category}. The word \textit{semi} means that the monoidal structure does not necessarily have a unit. This section contains also various calculations which will be used in the sequel. 
\item Section~\ref{basic-property-moore-flow} introduces the notion of \textit{$\mathcal{P}$-flow} and proves various basic properties, like its local presentability. It is also proved that the \textit{(execution) path $\mathcal{P}$-space functor} from $\mathcal{P}$-flows to $\mathcal{P}$-spaces is a right adjoint. 
\item Section~\ref{Isaev} recalls some results about Isaev's work \cite{Isaev} about model categories of fibrant objects. It also provides, as an easy consequence of Isaev's paper, a characterization of the class of weak equivalences as a small injectivity class with an explicit description of the generating set. The latter result is a particular case of a more general result due to Bourke.
\item Section~\ref{homotopytheoryMoore} expounds and describes as completely as possible the \textit{q-model structure} of $\mathcal{P}$-flows. It is very similar to the q-model structure of flows in many ways. In fact, we follow the method of \cite{leftdetflow} step by step constructing the q-model structure of $\dtop$ by using \cite{Isaev}.
\item Section~\ref{cofibrant-moore-flow} proves that the path $\mathcal{P}$-space functor from $\mathcal{P}$-flows to $\mathcal{P}$-spaces takes (trivial resp.) q-cofibrations of $\mathcal{P}$-flows between q-cofibrant $\mathcal{P}$-flows to (trivial resp.) projective q-cofibrations of $\mathcal{P}$-spaces. In particular, the $\mathcal{P}$-space of execution paths of a q-cofibrant $\mathcal{P}$-flow is a projective q-cofibrant $\mathcal{P}$-space. The material expounded in \cite{leftproperflow} is used in a crucial way. 
\item  Finally, Section~\ref{equivalencePflow-flow}, after a reminder about flows, establishes the Quillen equivalence between the q-model structure of $\mathcal{P}$-flows and the q-model structure of flows. The proof requires to use the main theorem of \cite{dgrtop} which is recalled in Theorem~\ref{eq-topdgr-top}. 
\end{itemize}

\subsection*{Notations, conventions and prerequisites}

We refer to \cite{TheBook} for locally presentable categories, to \cite{MR2506258} for combinatorial model categories.  We refer to \cite{MR99h:55031,ref_model2} for more general model categories. We refer to \cite{KellyEnriched} and to \cite[Chapter~6]{Borceux2} for enriched categories. All enriched categories are topologically enriched categories: \textit{the word topologically is therefore omitted}. What follows is some notations and conventions.

\begin{itemize}[leftmargin=*]
	\item $A:=B$ means that $B$ is the definition of $A$.
	\item $\K^{op}$ denotes the opposite category of $\K$.
	\item $\Obj(\K)$ is the class of objects of $\K$, $\Mor(\K)$ is the category of morphisms of $\K$ with the commutative squares for the morphisms.
	\item $\K^I$ is the category of functors and natural transformations from a small category $I$ to $\K$.
	\item $\varnothing$ is the initial object, $\mathbf{1}$ is the final object, $\id_X$ is the identity of $X$.
	\item $\K(X,Y)$ is the set of maps in a set-enriched, i.e. locally small, category $\K$.
	\item $\K(X,Y)$ is the space of maps in an enriched category $\K$. The underlying set of maps may be denoted by $\K_0(X,Y)$ if it is necessary to specify that we are considering the underlying set. 
	\item The composition of two maps $f:A\to B$ and $g:B\to C$ is denoted by $gf$ or, if it is helpful for the reader, by $g.f$; the composition of two functors is denoted in the same way.
\end{itemize}

\subsection*{History, informal comments and acknowledgments}

The route which would lead to this work is long. After the redaction of \cite{mdtop}, I thought that the intermediate category to put between multipointed $d$-spaces and flows was a Moore version of multipointed $d$-spaces (see e.g \cite{135738}). The latter idea seems to be a dead-end. In fall 2017, I had the idea to use a Moore version of flows instead of multipointed $d$-spaces. It did not come to my mind at first because the composition of paths is already strictly associative on flows. This led me to a first attempt to define Moore flows in spring 2018. A new obstacle then arised. Due to the non-contractibility of the classifying space of the group of nondecreasing homeomorphisms from $[0,1]$ to itself, the behavior of the homotopy colimit of spaces (see the introduction of \cite{dgrtop}) prevents the left Quillen functor $\lmoore$ from Moore flows to flows of Proposition~\ref{def-N} from being homotopically surjective. I had fruitful email discussions with Tim Porter during summer 2018, and I thank him for that, from which a better idea emerged: I had to work in an enriched setting. The preceding obstable is then overcome thanks to \cite[Theorem~7.6]{dgrtop} which is recalled in Theorem~\ref{eq-topdgr-top}. It was then necessary to prove that a q-cofibrant Moore flow has a projective q-cofibrant $\mathcal{P}$-space of execution paths. This fact is required also for the functor $\lmoore$ to be homotopically surjective. The redaction of this proof (see Theorem~\ref{PathSpacePreservesCofibrancy}) led me to the formulation of a Moore flow as a semicategory enriched over the biclosed semimonoidal structure of $\mathcal{P}$-spaces. Some explanations about the axioms satisfied by a general reparametrization category are given after Definition~\ref{def-reparam}.

I thank the two anonymous referees for their reports and suggestions. One of the two referees suggested two terminologies $\mathcal{P}$-space or $\mathcal{P}$-trace for Definition~\ref{pspace}. When $\mathcal{P}$ is the terminal category $\mathbf{1}$, the semimonoidal category $(\topdgr_0,\ot)$ is the (semi)monoidal category $(\top,\p)$. It is the reason why the first suggestion is chosen.

\section{Reminder about topological spaces}
\label{reminder-space}

The category $\top$ is either the category of \textit{$\Delta$-generated spaces} or the full subcategory of \textit{$\Delta$-Hausdorff $\Delta$-generated spaces} (cf. \cite[Section~2 and Appendix~B]{leftproperflow}). We summarize some basic properties of $\top$ for the convenience of the reader: 
\begin{itemize}[leftmargin=*]
	\item $\top$ is locally presentable.
	\item The inclusion functor from the full subcategory of $\Delta$-generated spaces to the category of general topological spaces together with the continuous maps has a right adjoint called the $\Delta$-kelleyfication functor. The latter functor does not change the underlying set.
	\item Let $A\subset B$ be a subset of a space $B$ of $\top$. Then $A$ equipped with the $\Delta$-kelleyfication of the relative topology belongs to $\top$. 
	\item The colimit in $\top$ is given by the final topology in the following situations: 
	\begin{itemize}
		\item A transfinite compositions of one-to-one maps.
		\item A pushout along a closed inclusion.
		\item A quotient by a closed subset or by an equivalence relation having a closed graph.
	\end{itemize}
	In these cases, the underlying set of the colimit is therefore the colimit of the underlying sets. In particular, the CW-complexes, and more generally all cellular spaces are equipped with the final topology. Note that cellular spaces are even Hausdorff (and paracompact, normal, etc...).
	\item $\top$ is cartesian closed. The internal hom $\ttop(X,Y)$ is given by taking the $\Delta$-kelleyfication of the compact-open topology on the set $\mathcal{T\!O\!P}(X,Y)$ of all continuous maps from $X$ to $Y$. 
\end{itemize}
The q-model structure~\footnote{We use the terminology of \cite{ParamHomTtheory}.} is denoted by $\top_q$. It is enriched over itself using the binary product. The terminology of cartesian enrichment is used sometimes. We keep denoting the \textit{set} of continuous maps from $X$ to $Y$ by $\top(X,Y)$ and the \textit{space} of maps from $X$ to $Y$ by $\ttop(X,Y)$ like in our previous papers (and not by $\top_0(X,Y)$ and $\top(X,Y)$ respectively).

A topological space is connected if and only if it is path-connected and every topological space is homeomorphic to the disjoint sum of its nonempty path-connected components \cite[Proposition~2.8]{mdtop}. The space $\CC(Z)$ is the space of path-connected components of $Z$ equipped with the final topology with respect to the canonical map $Z\to \CC(Z)$, which turns out to be the discrete topology by \cite[Lemma~5.8]{leftproperflow}. 

Note that the paper \cite{dgrtop}, which is used several times in this work, is written in the category of $\Delta$-generated spaces; it is still valid in the category of $\Delta$-Hausdorff $\Delta$-generated spaces. This point is left as an exercise for the idle mathematician: see \cite[Section~2 and Appendix~B]{leftproperflow} for any help about the topology of $\Delta$-generated spaces. 

\begin{nota}
	Let $n\geq 1$. Denote by \[\mathbf{D}^n = \{b\in \mathbb{R}^n, |b| \leq 1\}\] the $n$-dimensional disk, and by \[\mathbf{S}^{n-1}= \{b\in \mathbb{R}^n, |b| = 1\}\] the $(n-1)$-dimensional sphere. By convention, let $\mathbf{D}^{0}=\{0\}$ and $\mathbf{S}^{-1}=\varnothing$. 
\end{nota}

\section{Reminder about enriched presheaves of topological spaces}
\label{presheaf-over-enriched}

\begin{nota}
	$\mathcal{P}$ denotes an enriched small category in this section. The enriched category of enriched presheaves from $\mathcal{P}^{op}$ to $\top$ is denoted by $\topdgr$. The underlying set-enriched category of enriched maps of enriched presheaves is denoted by $\topdgr_0$. 
\end{nota}

The enriched category $\topdgr$ is tensored and cotensored \cite[Lemma~5.2]{MoserLyne}. For an enriched presheaf $F:\mathcal{P}^{op}\to \top$, and a topological space $U$, the enriched presheaf $F\ot U:\mathcal{P}^{op}\to \top$ is defined by $F\ot U = F(-)\p U$ and the enriched presheaf $F^U:\mathcal{P}^{op}\to \top$ is defined by $F^U=\ttop(U,F(-))$. By \cite[Definition~6.2.4 and Diagram~6.13]{Borceux2} (see \cite[Section~5]{dgrtop} for more detailed explanations), the underlying category $\topdgr_0$ of $\topdgr$ can be identified with the full subcategory of the category  $\top^{\mathcal{P}^{op}}$ of functors $F:\mathcal{P}^{op}\to \top$ such that the set map  $\mathcal{P}_0(\ell_1,\ell_2) \longrightarrow \top(F(\ell_2),F(\ell_1))$ induces a continuous map  $\mathcal{P}(\ell_1,\ell_2) \longrightarrow \ttop(F(\ell_2),F(\ell_1))$ for all $\ell_1,\ell_2\in \Obj(\mathcal{P})$. By \cite[Proposition~5.1]{dgrtop}, the category $\topdgr_0$ is locally presentable. The category $\topdgr_0$ is a full reflective and coreflective subcategory of $\top^{\mathcal{P}^{op}_0}$ (see e.g. \cite[Proposition~5.3]{dgrtop}). The enriched structure of $\topdgr$ combined with the enriched Yoneda lemma \cite[Theorem~6.3.5]{Borceux2} implies the natural homeomorphism 
\begin{equation}
\topdgr(\mathcal{P}(-,\ell)\p U,D) \iso \ttop(U,D(\ell)) \label{ey} \tag{En-Yo}
\end{equation}
for all $\ell\in \Obj(\mathcal{P})$, for all topological spaces $U$ and all enriched presheaves $D$. For any pair of enriched presheaves $D,E\in \topdgr$, by \cite[Equation~2.10]{KellyEnriched}, there is the homeomorphism
\begin{equation}
\topdgr(D,E) \iso \int_\ell \ttop(D(\ell),E(\ell)) \label{en} \tag{En-Nat}
\end{equation}
and consequently, together with the enriched Yoneda lemma, we obtain the homeomorphism 
\begin{equation}
D= \int^{\ell} \mathcal{P}(-,\ell) \p D(\ell). \label{rep} \tag{En-Rep}
\end{equation}
Note that the coend $\int^\ell F(\ell,\ell)$ of an enriched functor $F:\mathcal{P}\p \mathcal{P}^{op}\to \top$ is defined as the coequalizer of the two maps \cite[dual of Equation~2.2]{KellyEnriched}: 
\begin{equation}\label{coendeq}
\bigsqcup_{(\ell_1,\ell_2)} \mathcal{P}(\ell_1,\ell_2) \p F(\ell_1,\ell_2) \rightrightarrows \bigsqcup_\ell F(\ell,\ell).  \tag{Coend}
\end{equation}
It implies that for any enriched functor $F:\mathcal{P} \p \mathcal{P}^{op}\to \top$ and any topological space $U$, there is a homeomorphism 
\begin{equation}\label{cceq}
\bigg(\int^\ell F(\ell,\ell)\bigg)\p U \iso \int^\ell \big(F(\ell,\ell)\p U\big)  \tag{CC}
\end{equation}
because $\top$ is cartesian closed. 

\begin{nota}
	Let $U$ be a topological space. The constant presheaf $U$ is denoted by $\Delta_{\mathcal{P}^{op}}U$. It is enriched since for all $\ell_1,\ell_2\in \Obj(\mathcal{P})$ the map $\mathcal{P}(\ell_1,\ell_2)\to (\Delta_{\mathcal{P}^{op}}U)(\ell_1,\ell_2)$ is the constant map which is therefore continuous. The enriched presheaf $\Delta_{\mathcal{P}^{op}}\varnothing$ is called the {\rm empty enriched presheaf}.
\end{nota}

\begin{nota}
	$\liminj : \topdgr_0 \to \top$ denotes the colimit functor: there is no reason to specify in the notation the underlying category $\mathcal{P}^{op}$.
\end{nota}

Finally, by using the fact that the functor $\Delta_{\mathcal{P}^{op}}:\top\to \topdgr_0$ preserves cotensors, we deduce using \cite[Theorem~6.7.6]{Borceux2} that for every enriched presheaf $D$ and every topological space $Z$, there is the natural homeomorphism 
\begin{equation}\label{colim-adj}
\ttop(\liminj D,Z) \iso \topdgr(D,\Delta_{\mathcal{P}^{op}} Z).  \tag{En-Adj}
\end{equation}

\section{Reparametrization category}
\label{param}

\bd \label{def0} A {\rm semimonoidal category} $(\K,\ot)$ is a category $\K$ equipped with a functor $\ot:\K\p \K\to \K$ together with a natural isomorphism $a_{x,y,z}: (x\ot y) \ot z \to x \ot (y\ot z)$ called the {\rm associator} satisfying the pentagon axiom \cite[diagram~(5) page 158]{MR1712872}. 
\ed

According to the usual terminology used for similar situations, a semimonoidal category could be called a \textit{non-unital monoidal category}. Note that it is the monoidal structure which is non unital, not the category. Non unital categories (a.k.a semicategories) appear on stage in Section~\ref{basic-property-moore-flow} of this paper.

\bd \label{def1} A semimonoidal category $(\K,\ot)$ is {\rm enriched}~\footnote{Remember that in this paper, all enriched categories are enriched over $\top$.} if the category $\K$ is enriched and if the set map \[\K(a,b)\p \K(c,d) \longrightarrow \K(a\ot c,b\ot d)\] is continuous for all objects $a,b,c,d\in \Obj(\K)$. 
\ed

\bd \label{def-reparam}
A {\rm reparametrization category} $(\mathcal{P},\ot)$ is a small enriched semimonoidal category satisfying the following additional properties: 
\begin{enumerate}
	\item The semimonoidal structure is strict, i.e. the associator is the identity.
	\item All spaces of maps $\mathcal{P}(\ell,\ell')$ for all objects $\ell$ and $\ell'$ of $\mathcal{P}$ are contractible. 
	\item For all maps $\phi:\ell\to \ell'$ of $\mathcal{P}$, for all $\ell'_1,\ell'_2\in \Obj(\mathcal{P})$ such that $\ell'_1\ot\ell'_2=\ell'$, there exist two maps $\phi_1:\ell_1\to \ell'_1$ and $\phi_2:\ell_2\to \ell'_2$ of $\mathcal{P}$ such that $\phi=\phi_1 \ot \phi_2 : \ell_1\ot\ell_2 \to \ell'_1 \ot\ell'_2$ (which implies that $\ell_1 \ot \ell_2=\ell$). 
\end{enumerate}
\ed

The semimonoidal structure enables us to have a semigroup structure on objects, to formulate the third axiom and to define the enriched functors $s^L_\ell$ and $s^R_\ell$ of Proposition~\ref{dec}. The semigroup structure on objects is important for Lemma~\ref{petitcalcul}. It plays a central role for the study of the biclosed semimonoidal structure of $\mathcal{P}$-spaces defined in Definition~\ref{semimono} (see Theorem~\ref{closedsemi}). This biclosed semimonoidal structure is required to formalize $\mathcal{P}$-flows as enriched semicategories in Definition~\ref{PF}. This point of view on $\mathcal{P}$-flows, which is not the initial one I had chronologically, is used mainly for the proof of the key fact that a q-cofibrant $\mathcal{P}$-flow has a projective q-cofibrant $\mathcal{P}$-space of execution paths (Theorem~\ref{pretower} and Theorem~\ref{PathSpacePreservesCofibrancy}). The latter proof relies on the calculations made in \cite{leftproperflow} in the setting of semicategories enriched over topological spaces (a.k.a. flows). The reparametrization category must be enriched to be able to take into account the contractibility of the spaces of maps. All spaces of maps must be contractible to use Theorem~\ref{eq-topdgr-top}. Otherwise the non-contractibility of the classifying space of the group of nondecreasing homeomorphisms from $[0,1]$ to itself prevents the left Quillen functor $\lmoore$ of Proposition~\ref{def-N} from $\mathcal{P}$-flows to flows from being homotopically surjective. The third axiom is used in the proof of Proposition~\ref{Ptenseur}. It enables us in Section~\ref{equivalencePflow-flow} to define the right Quillen functor $\moore$ from flows to $\mathcal{P}$-flows in Proposition~\ref{def-preN}.

\begin{nota}
	To stick to the intuition, we set $\ell+\ell' := \ell \ot \ell'$ for all $\ell,\ell'\in \Obj(\mathcal{P})$. Indeed, morally speaking, $\ell$ is the length of a path.
\end{nota}

A reparametrization category $\mathcal{P}$ is an enriched category with contractible spaces of morphisms such that  the set $\Obj(\mathcal{P})$ of objects of $\mathcal{P}$ has a structure of a semigroup with a composition law denoted by $+$, such that the set map 
\[
\ot : \mathcal{P}(\ell_1,\ell'_1) \p \mathcal{P}(\ell_2,\ell'_2) \to \mathcal{P}(\ell_1+\ell_2,\ell'_1+\ell'_2)
\]
is continuous for all $\ell_1,\ell'_1,\ell_2,\ell'_2 \in \Obj(\mathcal{P})$, and such that every map of $\mathcal{P}$ is of the form $\phi_1 \ot \phi_2$ (not necessarily in a unique way).

\begin{exa} \label{terminal-reparam}
	The terminal category with one object $\underline{1}$ and one map $\id_{\underline{1}}$ is a reparametri\-zation category. 
\end{exa}

\bd Let $\phi_i:[0,\ell_i] \to [0,\ell'_i]$ for $i=1,2$ with $\ell_1,\ell'_1,\ell_2,\ell'_2>0$ be two continuous maps preserving the extrema where a notation like $[0,\ell]$ means a segment of the real line. Then the map
\[
\phi_1 \ot \phi_2 : [0,\ell_1+\ell_2] \to [0,\ell'_1 + \ell'_2]
\]
denotes the continuous map defined by 
\[
(\phi_1 \ot \phi_2)(t) = \begin{cases}
\phi_1(t) & \hbox{if } 0\leq t\leq \ell_1\\
\phi_2(t-\ell_1)+\ell'_1 & \hbox{if } \ell_1\leq t\leq \ell_1+\ell_2\\
\end{cases}
\] 
\ed

We have the obvious proposition:

\bp \label{assoc-tenseur}
Let $\phi_i:[0,\ell_i] \to [0,\ell'_i]$ for $1\leq i \leq 3$ preserving the extrema. Then 
\[
(\phi_1 \ot \phi_2) \ot \phi_3 = \phi_1 \ot (\phi_2 \ot \phi_3).
\]
\ep

\begin{nota}
	The notation $[0,\ell_1]\iso^+ [0,\ell_2]$ for two real numbers $\ell_1,\ell_2>0$ means a nondecreasing homeomorphism from $[0,\ell_1]$ to $[0,\ell_2]$. It takes $0$ to $0$ and $\ell_1$ to $\ell_2$. 
\end{nota}

\bp \label{paramG}
There exists a reparametrization category, denoted by ${\mathcal{G}}$, such that the semigroup of objects is the open interval $]0,+\infty[$ equipped with the addition and such that for every $\ell_1,\ell_2>0$, there is the equality \[\mathcal{G}(\ell_1,\ell_2)=\{[0,\ell_1]\iso^+ [0,\ell_2]\}\] where the topology is the $\Delta$-kelleyfication of the relative topology induced by the set inclusion $\mathcal{G}(\ell_1,\ell_2) \subset \ttop([0,\ell_1],[0,\ell_2])$ and such that for every $\ell_1,\ell_2,\ell_3>0$, the composition map \[\mathcal{G}(\ell_1,\ell_2)\p \mathcal{G}(\ell_2,\ell_3) \to \mathcal{G}(\ell_1,\ell_3)\] is induced by the composition of continuous maps.
\ep

\bpf
Let $\phi_i:[0,\ell_i] \iso^+ [0,\ell'_i]$ for $i=1,2$ be two continuous maps. Then the map $\phi_1 \ot \phi_2 : [0,\ell_1+\ell_2] \to [0,\ell'_1 + \ell'_2]$ is a nondecreasing homeomorphism. The set map 
\[
\ot : \mathcal{G}(\ell_1,\ell'_1) \p \mathcal{G}(\ell_2,\ell'_2) \to \mathcal{G}(\ell_1+\ell_2,\ell'_1+\ell'_2)
\]
is continuous because the set map 
\[
[0,\ell_1+\ell_2] \p \mathcal{G}(\ell_1,\ell'_1) \p \mathcal{G}(\ell_2,\ell'_2) \to [0,\ell'_1+\ell'_2]
\]
takes $(t,\phi_1,\phi_2)$ to $(\phi_1\ot \phi_2)(t)$ which is continuous. The pentagon axiom is clearly satisfied. For every $\ell_1,\ell_2>0$, the continuous mapping $(t,\phi,u)\mapsto (1-t)\phi(u) + t(u\ell_2/\ell_1)$ yields a homotopy from $\id_{\mathcal{G}(\ell_1,\ell_2)}$ to the constant map taking any map of $\mathcal{G}(\ell_1,\ell_2)$ to the homothetie $u\mapsto u\ell_2/\ell_1$. Therefore the space $\mathcal{G}(\ell_1,\ell_2)$ is contractible. Let $\phi:[0,\ell]\iso^+ [0,\ell']$. Let $\ell'=\ell'_1+\ell'_2$. Let $\ell_1=\phi^{-1}(\ell'_1)$ and $\ell_2=\ell-\ell_1$. Define 
\[
\begin{cases}
\phi_1(t) = \phi(t) & \hbox{for $t\in [0,\ell_1]$}\\
\phi_2(t) = \phi(\ell_1+t) - \ell'_1 & \hbox{for $t\in [0,\ell_2]$}\\
\end{cases}
\]
Then $\phi_1:[0,\ell_1]\iso^+ [0,\ell'_1]$, $\phi_2:[0,\ell_2]\iso^+ [0,\ell'_2]$ and $\phi=\phi_1 \ot \phi_2$.
\epf

\begin{nota}
	The notation $[0,\ell_1]\to^+ [0,\ell_2]$ for two real numbers $\ell_1,\ell_2>0$ means a nondecreasing continuous from $[0,\ell_1]$ to $[0,\ell_2]$ preserving the extrema. It takes $0$ to $0$ and $\ell_1$ to $\ell_2$. 
\end{nota}

\bp \label{paramM}
There exists a reparametrization category, denoted by ${\mathcal{M}}$, such that the semigroup of objects is the open interval $]0,+\infty[$ equipped with the addition and such that for every $\ell_1,\ell_2>0$, there is the equality \[\mathcal{M}(\ell_1,\ell_2)=\{[0,\ell_1]\longrightarrow^+ [0,\ell_2]\}\] where the topology is the $\Delta$-kelleyfication of the relative topology induced by the set inclusion $\mathcal{M}(\ell_1,\ell_2) \subset \ttop([0,\ell_1],[0,\ell_2])$ and such that for every $\ell_1,\ell_2,\ell_3>0$, the composition map \[\mathcal{M}(\ell_1,\ell_2)\p \mathcal{M}(\ell_2,\ell_3) \to \mathcal{M}(\ell_1,\ell_3)\] is induced by the composition of continuous maps.
\ep

\bpf
The proof is similar to the proof of Proposition~\ref{paramG}.
\epf

In the cases of $(\mathcal{G},+)$ and $(\mathcal{M},+)$, the functors $(\ell,\ell')\mapsto \ell+\ell'$ and $(\ell,\ell') \mapsto \ell'+\ell$ coincide on objects, but not on morphisms. The terminal category is a symmetric reparametrization category. We do not know if there exist symmetric reparametrization categories not equivalent to the terminal category.

\section{\texorpdfstring{$\mathcal{P}$}{Lg}-spaces}
\label{presheaf-over-param}

From now on, $\mathcal{P}$ denotes a reparametrization category.

\bdn \label{pspace} An object of $\topdgr_0$ is called a {\rm $\mathcal{P}$-space}. Let $D$ be a $\mathcal{P}$-space. Let $\phi:\ell\to \ell'$ be a map of $\mathcal{P}$. Let $x\in D(\ell')$. We will use the notation 
	\[
	x.\phi := D(\phi)(x).
	\]
	The motivating example is if $x$ is a path of length $\ell'$, then $x.\phi$ is a path of length $\ell$ which is the reparametrization by $\phi$ of $x$.
\edn

\bd A $\mathcal{P}$-space $D$ is {\rm path-connected} if for every $\ell\in \mathcal{P}$, $D(\ell)$ is path-connected.
\ed

\begin{lem} \label{CCproduct}
The connected component functor $\CC:\top \to \set$, where $\set$ is the category of sets, preserves binary products.
\end{lem}

\bpf 
Let $X$ and $Y$ be two topological spaces. There is a canonical set map 
\[
\CC(X\p Y) \longrightarrow \CC(X) \p \CC(Y)
\] 
induced by the universal property of the binary product. This set map is clearly onto. Let $(x_0,y_0)$ and $(x_1,y_1)$ be two points of $X\p Y$ such that $x_0$ and $x_1$ ($y_0$ and $y_1$ resp.) are in the same connected components. Then there exist  continuous paths $\phi:[0,1]\to X$ and $\psi:[0,1]\to Y$ such that $\phi(i)=x_i$ and $\psi(i)=y_i$ for $i=0,1$. Therefore $(\phi,\psi)$ is a continuous path from $(x_0,y_0)$ to $(x_1,y_1)$. Thus $(x_0,y_0)$ and $(x_1,y_1)$ are in the same connected components. 
\epf

\begin{nota}
	Let $\mathcal{K}$ be an enriched category. Denote by $\thin{\mathcal{K}}$ the enriched category with the same objects as $\mathcal{K}$ and with \[\thin{\mathcal{K}}(\ell,\ell') = \CC(\mathcal{K}(\ell,\ell'))\] for all $\ell,\ell'\in \Obj(\mathcal{K})$. The composition law is defined thanks to Lemma~\ref{CCproduct} by the composite map: 
	\[
	\thin{\mathcal{K}}(\ell,\ell') \p \thin{\mathcal{K}}(\ell',\ell'') \iso \CC(\mathcal{K}(\ell,\ell')\p \mathcal{K}(\ell',\ell'')) \longrightarrow \thin{\mathcal{K}}(\ell,\ell'').
	\]
\end{nota}

\bp \label{cc-dgrtop}
Every $\mathcal{P}$-space $D$ can be decomposed as a coproduct
\[
D = \bigsqcup_{c} D_c
\]
with $D_c \neq \Delta_{\mathcal{P}}\varnothing$ and path-connected.
\ep

\bpf
For every map $\ell,\ell'\in \Obj(\mathcal{P})$ and every $\mathcal{P}$-space $D$ such that $D\neq \Delta_{\mathcal{P}}\varnothing$, the continuous map \[\mathcal{P}(\ell,\ell') \p D(\ell')\longrightarrow D(\ell)\] induces by Lemma~\ref{CCproduct} a composite continuous map 
\[
\mathcal{P}(\ell,\ell') \p \CC(D(\ell')) \to \CC(\mathcal{P}(\ell,\ell')) \p \CC(D(\ell')) \iso \CC(\mathcal{P}(\ell,\ell') \p D(\ell')) \to \CC(D(\ell)).
\]
Since $\top$ is cartesian closed, we obtain a continuous map 
\[
\mathcal{P}(\ell,\ell') \longrightarrow \ttop (\CC(D(\ell')),\CC(D(\ell)))
\]
from the contractible space $\mathcal{P}(\ell,\ell')$ to the discrete space $\ttop (\CC(D(\ell')),\CC(D(\ell)))$ which is therefore constant. Consequently, the set map \[\CC(D(\phi)):\CC(D(\ell'))\longrightarrow \CC(D(\ell))\] does not actually depend of $\phi \in \mathcal{P}(\ell,\ell')$. We obtain the commutative diagram of functors 
\[
\xymatrix@C=3em@R=3em
{
	\mathcal{P}^{op} \fr{D} \fd{}&\top \ar@{->}[d]^-{\CC} \\
	\thin{\mathcal{P}^{op}} \ar@{-->}[r]^-{\exists ! } & \set
}
\]
It turns out that the category $\thin{\mathcal{P}^{op}}$ is equivalent to the terminal category $\mathbf{1}$ with one object $\underline{1}$ and one map $\id_{\underline{1}}$. We obtain a commutative diagram of functors 
\[
\xymatrix@C=3em@R=3em
{
	\mathcal{P}^{op} \fr{D} \fd{}&\top \ar@{->}[d]^-{\CC} \\
	\mathbf{1} \ar@{-->}[r]^-{\exists ! \overline{D}} & \set
}
\]
Each element of $\overline{D}(\underline{1})$ corresponds for every object $\ell$ of $\mathcal{P}$ to a path-connected component of $D(\ell)$. The above commutative diagram also tells that for every map $\phi:\ell\to \ell'$ of $\mathcal{P}$, the map \[\CC(D(\phi)):\CC(D(\ell'))\longrightarrow \CC(D(\ell))\] is the identity of $\overline{D}(\underline{1})$. Thus there exists a decomposition 
\[
D = \bigsqcup_{c\in \overline{D}(\underline{1})} D_c
\]
in $\topdgr_0$ such that for every $\ell\in \Obj(\mathcal{P})$, the space $D_c(\ell)$ is nonempty path-connected.
\epf

\bd 
The set \[\CC(D):=\{D_c\mid c\in\overline{D}(\underline{1})\}\] is called the set of {\rm (nonempty) path-connected components} of $D$. By convention, \[\CC(\Delta_{\mathcal{P}^{op}}\varnothing)=\varnothing.\]
By considering the particular case where the reparametrization category $\mathcal{P}$ is the terminal category, we see that this generalizes the notion of path-connected component of a topological space. 
\ed

\bd \label{semimono}
Let $D$ and $E$ be two $\mathcal{P}$-spaces. Let 
\[
D \ot E = \int^{(\ell_1,\ell_2)} \mathcal{P}(-,\ell_1+\ell_2) \p D(\ell_1) \p E(\ell_2).
\]
\ed

\bp \label{dec} (The left shift functor and the right shift functor)
The following data assemble to an enriched functor $s^L_\ell:\mathcal{P}\to \mathcal{P}$: 
\[
\begin{cases}
s^L_\ell(\ell')= \ell + \ell'\\
s^L_\ell(\phi) = \id_\ell \ot \phi & \hbox{for a map $\phi:\ell'\to \ell''$.}
\end{cases}
\]
The following data assemble to an enriched functor $s^R_\ell:\mathcal{P}\to \mathcal{P}$: 
\[
\begin{cases}
s^R_\ell(\ell')= \ell' + \ell\\
s^R_\ell(\phi) =  \phi \ot \id_\ell & \hbox{for a map $\phi:\ell'\to \ell''$.}
\end{cases}
\]
\ep

\bpf
The maps $\mathcal{P}(\ell',\ell'') \to \mathcal{P}(\ell+\ell',\ell+\ell'')$ and $\mathcal{P}(\ell',\ell'') \to \mathcal{P}(\ell'+\ell,\ell''+\ell)$ are continuous for all $\ell',\ell''\in \Obj(\mathcal{P})$. 
\epf


\begin{nota} \label{nota0}
	With $D\in \topdgr_0$, set $(s^L_\ell)^*(D) = D.s^L_\ell \in \topdgr_0$ and $(s^R_\ell)^*(D) = D.s^R_\ell \in \topdgr_0$.
\end{nota}

%

\begin{lem} \label{petitcalcul} 
	For all $\ell',\ell''\in \Obj(\mathcal{P})$, there are the isomorphisms of $\mathcal{P}$-spaces (natural with respect to $\ell'$ and $\ell''$)
	\begin{align*}
	&&\int^{\ell} \mathcal{P}(-,\ell+\ell') \p \mathcal{P}(\ell,\ell'') \iso \mathcal{P}(-,\ell''+\ell'),\\
	&&\int^{\ell} \mathcal{P}(-,\ell'+\ell) \p \mathcal{P}(\ell,\ell'') \iso \mathcal{P}(-,\ell'+\ell'').
	\end{align*}
	The first isomorphism takes the equivalence class of $(\psi,\phi)\in \mathcal{P}(-,\ell+\ell') \p \mathcal{P}(\ell,\ell'')$ to $(s^R_{\ell'})^*(\phi)\psi = (\phi\ot \id_{\ell'})\psi$. The second isomorphism takes the equivalence class of $(\psi,\phi) \in \mathcal{P}(-,\ell'+\ell) \p \mathcal{P}(\ell,\ell'')$ to $(s^L_{\ell'})^*(\phi)\psi =(\id_{\ell'}\ot \phi)\psi$.
\end{lem}

\bpf
Pick a $\mathcal{P}$-space $D$. Then there is the sequence of homeomorphisms
\[\begin{aligned}
\topdgr\bigg(\int^{\ell} \mathcal{P}(-,\ell+\ell') \p \mathcal{P}(\ell,\ell''),D\bigg)
&\iso \int_{\ell} \topdgr\big(\mathcal{P}(-,\ell+\ell') \p \mathcal{P}(\ell,\ell''),D\big)\\
&\iso \int_{\ell} \ttop(\mathcal{P}(\ell,\ell''),D(\ell+\ell'))\\
&\iso \topdgr(\mathcal{P}(-,\ell''),(s^R_{\ell'})^*D)\\
&\iso D(\ell''+\ell')\\
&\iso \topdgr(\mathcal{P}(-,\ell''+\ell'),D),
\end{aligned}\]
the first homeomorphism since $\topdgr$ is enriched, the second the fourth and the fifth homeomorphisms by \eqref{ey} and the third homeomorphism by \eqref{en}. By composing with the functor $\top(\{0\},-)$, we obtain the natural bijection of sets 
\[
\topdgr_0\bigg(\int^{\ell} \mathcal{P}(-,\ell+\ell') \p \mathcal{P}(\ell,\ell''),D\bigg)\iso \topdgr_0(\mathcal{P}(-,\ell''+\ell'),D).
\]
The proof of the first isomorphism is complete thanks to the Yoneda lemma. The proof of the second isomorphism is similar and is left to the reader.
\epf

\bp \label{asso}
Let $D_1$ and $D_2$ be two $\mathcal{P}$-spaces and $L\in\Obj(\mathcal{P})$. Then the mapping $(x,y) \mapsto (\id,x,y)$ yields a surjective continuous map 
\[
\displaystyle\bigsqcup\limits_{\substack{(\ell_1,\ell_2)\\\ell_1+\ell_2=L}}  D_1(\ell_1)\p  D_2(\ell_2)\longrightarrow (D_1\ot D_2)(L).
\]
Moreover, the functor \[\ot : \topdgr_0\p \topdgr_0 \to \topdgr_0\] induces a semimonoidal structure on $\topdgr_0$. 
\ep

\bpf 
Let $(\psi,x_1,x_2) \in \mathcal{P}(L,\ell_1+\ell_2) \p D_1(\ell_1)\p  D_2(\ell_2)$ be a representative of an element of $(D_1\ot D_2)(L)$. Then there exist two maps $\psi_i:\ell'_i\to \ell_i$ for $i=1,2$ such that $\psi=\psi_1\ot \psi_2$. By definition of a coend (see Corollary~\ref{expl-descr} for a more detailed explanation), one has $(\psi,x_1,x_2)\sim (\id_L,x_1\psi_1,x_2\psi_2)$ in $(D_1\ot D_2)(L)$ and the proof of the first statement is complete.

Let $D_1,D_2,D_3$ be three $\mathcal{P}$-spaces. Let $a_{D_1,D_2,D_3}:(D_1\ot D_2) \ot D_3 \to D_1 \ot (D_2\ot D_3)$ be the composite of the isomorphisms (by using Lemma~\ref{petitcalcul} twice)
\[\begin{aligned}
(D_1\ot &D_2) \ot D_3\\ &\iso  \int^{(\ell_1,\ell_2,\ell_3)} \bigg(\int^{\ell} \mathcal{P}(-,\ell+\ell_3) \p \mathcal{P}(\ell,\ell_1+\ell_2)\bigg)\p D_1(\ell_1) \p D_2(\ell_2) \p D_3(\ell_3)\\
&\iso \int^{(\ell_1,\ell_2,\ell_3)} \mathcal{P}(-,\ell_1+\ell_2+\ell_3) \p D_1(\ell_1) \p D_2(\ell_2) \p D_3(\ell_3)\\
&\iso \int^{(\ell_1,\ell_2,\ell_3)} \bigg(\int^{\ell} \mathcal{P}(-,\ell_1+\ell) \p \mathcal{P}(\ell,\ell_2+\ell_3)\bigg)\p D_1(\ell_1) \p D_2(\ell_2) \p D_3(\ell_3)\\
&\iso D_1 \ot (D_2\ot D_3).
\end{aligned}\]
Let $(\psi,(\phi,x_1,x_2),x_3)\in ((D\ot E)\ot F)(L)$ with $x_i\in D_i(\ell_i)$ for $i=1,2,3$ and $L\in\Obj(\mathcal{P})$. Write $\phi=\phi_1\ot\phi_2$ with $\phi_i:\ell'_i\to \ell_i$ for $i=1,2$ and $\psi=\psi_1\ot \psi_2 \ot \psi_3$ with $\psi_i:\ell''_i\to \ell'_i$ for $i=1,2,3$ with $\ell'_3=\ell_3$. In particular, $L=\ell''_1+\ell''_2+\ell''_3$. We obtain \[(\psi,(\phi,x_1,x_2),x_3) \sim (\id_{L},(\id_{\ell''_1+\ell''_2},x_1\phi_1\psi_1,x_2\phi_2\psi_2),x_3\psi_3)\] in $((D\ot E)\ot F)(L)$. The above sequence of isomorphisms takes the latter at first to the equivalence class of \[((\id_{\ell''_1+\ell''_2}\ot \id_{\ell''_3})\id_{L},x_1\phi_1\psi_1,x_2\phi_2\psi_2,x_3\psi_3)\] by Lemma~\ref{petitcalcul}, and, since $(\id_{\ell''_1+\ell''_2}\ot \id_{\ell''_3})\id_{L} = (\id_{\ell''_1}\ot \id_{\ell''_2 +\ell''_3})\id_{L}$ and by Lemma~\ref{petitcalcul} again, to the equivalence class of \[(\id_{L},x_1\phi_1\psi_1,(\id_{\ell''_2+\ell''_3},x_2\phi_2\psi_2,x_3\psi_3)).\] We deduce that the associator $a_{D,E,F}:(D\ot E) \ot F \to D \ot (E\ot F)$ satisfies the pentagon axiom thanks to the first part of the proof.
\epf

\bp \label{assoc-n}
Let $D_1,\dots,D_n$ be $n$ $\mathcal{P}$-spaces with $n\geq 1$. Then there is the natural isomorphism of $\mathcal{P}$-spaces
\[
D_1 \ot \dots \ot D_n \iso \int^{(\ell_1,\dots,\ell_n)} \mathcal{P}(-,\ell_1+\dots + \ell_n) \p D_1(\ell_1) \p \dots D_n(\ell_n).
\]
\ep

\bpf
Let us prove this isomorphism by induction for $n\geq 1$. The formula is satisfied for $n=1$ by \eqref{rep}. Suppose that it is proved for $n\geq 1$. Then there is the sequence of natural isomorphisms of $\mathcal{P}$-spaces 
\[\begin{aligned}
&(D_1 \ot \dots \ot D_n) \ot D_{n+1} \\
&\iso \int^{(\ell,\ell_{n+1})} \mathcal{P}(-,\ell+\ell_{n+1}) \p (D_1 \ot \dots \ot D_n)(\ell) \p D_{n+1}(\ell_{n+1})\\
&\iso \int^{(\ell_1,\dots,\ell_{n+1})} \bigg(\int^{\ell} \mathcal{P}(-,\ell + \ell_{n+1}) \p \mathcal{P}(\ell,\ell_1+\dots + \ell_n)\bigg) \p D_1(\ell_1) \p \dots \p  D_{n+1}(\ell_{n+1})\\
&\iso \int^{(\ell_1,\dots,\ell_{n+1})} \mathcal{P}(-,\ell_1+\dots + \ell_{n+1}) \p D_1(\ell_1) \p \dots \p  D_{n+1}(\ell_{n+1}),
\end{aligned}\]
the first and the second isomorphisms by definition of $\ot$ by Fubini and by \eqref{cceq}, and the third isomorphism by Lemma~\ref{petitcalcul}. The induction is complete.
\epf

\begin{cor} \label{expl-descr}
Let $D_1,\dots,D_n$ be $n$ $\mathcal{P}$-spaces with $n\geq 1$. Then for all $L\in \Obj(\mathcal{P})$, the space $(D_1 \ot \dots \ot D_n)(L)$ is the quotient of the space 
\[
\bigsqcup_{(\ell_1,\dots,\ell_n)} \mathcal{P}(L,\ell_1+\dots + \ell_n) \p D_1(\ell_1) \p \dots D_n(\ell_n)
\]
by the equivalence relation generated by the identifications 
\[
(\psi,x_1.\phi_1,\dots,x_n.\phi_n) = ((\phi_1\ot \dots \ot\phi_n).\psi,x_1,\dots,x_n)
\]
for all $\psi\in \mathcal{P}(L,\ell_1+\dots+\ell_n)$, all $x_i\in D_i(\ell'_i)$ and all $\phi_i\in \mathcal{P}(\ell_i,\ell'_i)$, where \[(\psi,x_1.\phi_1,\dots,x_n.\phi_n)\in \mathcal{P}(L,\ell_1+\dots + \ell_n) \p D_1(\ell_1) \p \dots D_n(\ell_n)\] and
\[((\phi_1\ot \dots \phi_n).\psi,x_1,\dots,x_n)\in \mathcal{P}(L,\ell'_1+\dots + \ell'_n) \p D_1(\ell'_1) \p \dots D_n(\ell'_n).\]
\end{cor}

\bpf Using \eqref{coendeq} and Proposition~\ref{assoc-n}, we see that the space $(D_1 \ot \dots \ot D_n)(L)$ is the coequalizer of the two maps
\[
\xymatrix@R=3em
	{
	\displaystyle\bigsqcup\limits_{\substack{(\ell_1,\dots,\ell_n)\\(\ell'_1,\dots,\ell'_n)}}
	\mathcal{P}(L,\ell_1+\dots + \ell_n)\p \big(\prod\limits_{1\leq i \leq n}\mathcal{P}(\ell_i,\ell'_i) \p D_i(\ell'_i)\big)   \ar@{->}@<-1ex>[d]_-{f_1} \ar@{->}@<1ex>[d]^-{f_2}\\
	\displaystyle\bigsqcup\limits_{(\ell_1,\dots,\ell_n)}\mathcal{P}(L,\ell_1+\dots + \ell_n) \p \big(\prod\limits_{1\leq i \leq n} D_i(\ell_i)\big)
	}
\]
defined by 
\[
f_1:(\psi,x_1,\phi_1,\dots,x_n,\phi_n) \mapsto ((\phi_1\ot \dots \ot\phi_n).\psi,x_1,\dots,x_n)
\]
and
\[
f_2:(\psi,x_1,\phi_1,\dots,x_n,\phi_n) \mapsto (\psi,x_1.\phi_1,\dots,x_n.\phi_n).
\]
Morally speaking, the $x_i's$ are execution paths. The map $f_1$ composes the reparametrizations $\phi_1,\dots,\phi_n$ with $\psi$ before applying them to the Moore composition $x_1*\dots*x_n$. The map $f_2$ reparametrizes each $x_i$ by $\phi_i$ and then reparametrizes the Moore composition $(x_1.\phi_1)*\dots*(x_n.\phi_n)$ by $\psi$.
\epf

When $\mathcal{P}$ is the terminal category with one object $\underline{1}$ and one map $\id_{\underline{1}}$, the tensor product $(D \ot E)(\underline{1})$ is the quotient of the space $D(\underline{1})\p E(\underline{1})$ by the discrete equivalence relation by Corollary~\ref{expl-descr}. In other terms, the tensor product coincides with the binary product and it has therefore a unit in this case. When $\mathcal{P}$ is the reparametrization category of Proposition~\ref{paramG} or of Proposition~\ref{paramM}, it is unlikely that the corresponding tensor product of $\mathcal{P}$-spaces over them has a unit but we are unable to prove it.

\bth \label{closedsemi} 
Let $D$, $E$ and $F$ be three $\mathcal{P}$-spaces. Let  
\begin{align*}
&\{E,F\}_L:= \ell\mapsto \topdgr(E,(s^L_\ell)^*F),\\
&\{E,F\}_R:= \ell\mapsto \topdgr(E,(s^R_\ell)^*F).
\end{align*}
These yield two $\mathcal{P}$-spaces and there are the natural homeomorphisms
\begin{align*}
&\topdgr(D,\{E,F\}_L) \iso \topdgr(D\ot E,F),\\
&\topdgr(E,\{D,F\}_R) \iso \topdgr(D\ot E,F).
\end{align*}
Consequently, the functor \[\ot : \topdgr_0\p \topdgr_0 \to \topdgr_0\] induces a structure of biclosed semimonoidal structure on $\topdgr_0$.
\eth

\bpf
We only treat the case of $\{E,F\}_L$. The other case is similar. Since $(s^R_{L})^*F$ is an enriched functor for every $L\in \Obj(\mathcal{P})$ by Proposition~\ref{dec}, the set map \[\mathcal{P}(\ell,\ell') \longrightarrow \ttop(F(\ell'+L),F(\ell+L))\] is continuous for all $\ell,\ell'\in \Obj(\mathcal{P})$. Using the cartesian closedness of $\top$, we obtain a continuous map 
\[
\mathcal{P}(\ell,\ell') \p F(\ell'+L) \longrightarrow F(\ell+L)
\]
which is natural with respect to $L\in \Obj(\mathcal{P})$, i.e. an element $f_{\ell,\ell'}$ of the space of natural transformations
\[
\topdgr\big(\mathcal{P}(\ell,\ell') \p (s^L_{\ell'})^*F,(s^L_{\ell})^*F\big).
\]
The latter space is homeomorphic to 
\[
\ttop\big(\mathcal{P}(\ell,\ell'),\topdgr((s^L_{\ell'})^*F,(s^L_{\ell})^*F)\big)
\]
because $\topdgr$ is enriched. We obtain, by composition, a continuous map 
\[
\xymatrix@C=2.5em{\mathcal{P}(\ell,\ell') \p \{E,F\}_L(\ell') \fr{(f_{\ell,\ell'},\id)}& \topdgr((s^L_{\ell'})^*F,(s^L_{\ell})^*F) \p \{E,F\}_L(\ell') \fr{}& \{E,F\}_L(\ell),}
\]
where the right-hand map is induced by the composition of natural transformations (which is continuous). It means that the mapping $\ell\mapsto \topdgr(E,(s^L_\ell)^*F)$ yields a well-defined $\mathcal{P}$-space. This proves the first part of the statement of the theorem. There is the sequence of natural homeomorphisms 
\[\begin{aligned}
\topdgr(D,\{E,F\}_L) & \iso \int_\ell \ttop\big(D(\ell),\topdgr(E,(s^L_\ell)^*F)\big)\\
&\iso \int_{(\ell,\ell')} \ttop\big(D(\ell),\ttop(E(\ell'),F(\ell+\ell'))\big) \\
&\iso \int_{(\ell,\ell')} \ttop\big(D(\ell)\p E(\ell'),F(\ell+\ell')\big) \\
&\iso \int_{(\ell,\ell')} \topdgr\big(\mathcal{P}(-,\ell+\ell') \p D(\ell)\p E(\ell'),F\big)\\
&\iso \topdgr(D\ot E,F),
\end{aligned}\]
the first and second homeomorphisms by \eqref{en}, the third homeomorphism because $\top$ is enriched cartesian closed, the fourth homeomorphism by \eqref{ey}, and finally the last homeomorphism since $\topdgr$ is enriched and by definition of $\ot$. By composing with the functor $\top(\{0\},-)$, we obtain the desired adjunction. The proof is complete thanks to Proposition~\ref{asso}.
\epf

\begin{nota}
	Let \[\mathbb{F}^{\mathcal{P}^{op}}_{\ell}U=\mathcal{P}(-,\ell)\p U \in \topdgr_0\] where $U$ is a topological space and where $\ell$ is an object of $\mathcal{P}$.
\end{nota}

\bp \label{Ftenseur} 
Let $U,U'$ be two topological spaces. Let $\ell,\ell'\in \Obj(\mathcal{P})$. There is the natural isomorphism of $\mathcal{P}$-spaces 
\[
\mathbb{F}^{\mathcal{P}^{op}}_{\ell}U \ot \mathbb{F}^{\mathcal{P}^{op}}_{\ell'}U' \iso \mathbb{F}^{\mathcal{P}^{op}}_{\ell+\ell'}(U\p U').
\]
\ep

\bpf
One has 
\[
\mathbb{F}^{\mathcal{P}^{op}}_{\ell}U \ot \mathbb{F}^{\mathcal{P}^{op}}_{\ell'}U' = \int^{(\ell_1,\ell_2)} \mathcal{P}(-,\ell_1+\ell_2) \p \mathcal{P}(\ell_1,\ell) \p \mathcal{P}(\ell_2,\ell') \p U \p U'.
\]
Using Lemma~\ref{petitcalcul}, we obtain 
\[
\mathbb{F}^{\mathcal{P}^{op}}_{\ell}U \ot \mathbb{F}^{\mathcal{P}^{op}}_{\ell'}U' = \int^{\ell_1} \mathcal{P}(\ell_1,\ell) \p \mathcal{P}(-,\ell_1+\ell') \p U \p U'.
\]
Using Lemma~\ref{petitcalcul} again, we obtain
\[
\mathbb{F}^{\mathcal{P}^{op}}_{\ell}U \ot \mathbb{F}^{\mathcal{P}^{op}}_{\ell'}U' = \mathcal{P}(-,\ell+\ell') \p U \p U'.
\]
\epf

\bp \label{Ptenseur} 
Let $U$ and $U'$ be two topological spaces.  There is the natural isomorphism of $\mathcal{P}$-spaces
\[
\Delta_{\mathcal{P}^{op}} U \ot \Delta_{\mathcal{P}^{op}} U' \iso \Delta_{\mathcal{P}^{op}} (U \p U').
\]
\ep

\bpf
Since $\top$ is cartesian closed, it suffices to consider the case where $U=U'$ is a singleton. In that case, by Corollary~\ref{expl-descr}, the topological space $(\Delta_{\mathcal{P}^{op}} U \ot \Delta_{\mathcal{P}^{op}} U')(L)$ is the quotient of the space \[\bigsqcup_{(\ell,\ell')} \mathcal{P}(L,\ell+\ell')\] by the identifications $(\phi_1 \ot \phi_2).\phi \sim \phi$. Let $\psi\in \mathcal{P}(L,\ell+\ell')$ for some $\ell,\ell'\in \Obj(\mathcal{P})$. By definition of a reparametrization category, write $\psi=\psi_1 \ot \psi_2$ with $\psi_1:\ell_1\to \ell$ and $\psi_2:\ell_2\to \ell'$. Then we obtain $\psi = (\psi_1 \ot \psi_2).\id_L$. We deduce that $\psi\sim\id_L$ in $(\Delta_{\mathcal{P}^{op}} U \ot \Delta_{\mathcal{P}^{op}} U')(L)$. 
\epf

\bp \label{tensor-product} 
Let $D$ and $E$ be two $\mathcal{P}$-spaces. Then there is a natural homeomorphism 
\[
\liminj (D \ot E) \iso \liminj D \p \liminj E.
\]
\ep

\bpf
Let $Z$ be a topological space. There is the sequence of natural homeomorphisms 
\[\begin{aligned}
\ttop\big(\liminj (D \ot E),Z\big) & \iso \topdgr\big(D \ot E,\Delta_{\mathcal{P}^{op}}Z\big)\\
&\iso \topdgr \bigg(D,\ell \mapsto \topdgr(E,(s^L_\ell)^*\Delta_{\mathcal{P}^{op}}(Z))\bigg)\\
&\iso \topdgr \bigg(D,\Delta_{\mathcal{P}^{op}}\big(\topdgr(E,\Delta_{\mathcal{P}^{op}}(Z))\big)\bigg)\\
&\iso \ttop\big(\liminj D,\topdgr(E,\Delta_{\mathcal{P}^{op}}(Z))\big)\\
&\iso \ttop\big(\liminj D,\ttop(\liminj E,Z)\big)\\
&\iso \ttop\big((\liminj D) \p (\liminj E),Z \big),
\end{aligned}\]
the first fourth and fifth homeomorphisms by \eqref{colim-adj}, the second homeomorphism by Theorem~\ref{closedsemi}, the third homeomorphism since $(s^L_\ell)^*\Delta_{\mathcal{P}^{op}}(Z) = \Delta_{\mathcal{P}^{op}}(Z)$ for all $\ell \in \Obj(\mathcal{P})$ and the last homeomorphism since $\top$ is enriched cartesian closed. By composing with the functor $\top(\{0\},-)$, we obtain the natural bijection of sets 
\[
\top\big(\liminj (D \ot E),Z\big)\iso \top\big((\liminj D) \p (\liminj E),Z \big).
\]
The proof is complete thanks to the Yoneda lemma.
\epf

\section{\texorpdfstring{$\mathcal{P}$}{Lg}-flows}
\label{basic-property-moore-flow}

A \textit{semicategory}, also called \textit{nonunital category} in the literature, is a category without identity maps in the structure. It is \textit{enriched} over a biclosed monoidal category $(\mathcal{V},\ot,I)$ if it satisfies all axioms of enriched category except the one involving the identity maps, i.e. the enriched composition is associative and not necessarily unital. Since the existence of a unit $I$ is not necessary anymore, it makes sense to define the notion of semicategory enriched over a biclosed semimonoidal category. 

This section starts on purpose with the following very concise definition which is going to be explained right after.

\bd \label{PF} Let $\mathcal{P}$ be a reparametrization category. A {\rm $\mathcal{P}$-flow $X$} is a small semicategory enriched over the biclosed semimonoidal category $(\topdgr_0,\ot)$. The category of $\mathcal{P}$-flows is denoted by $\dtopP$. 
\ed

The following definition is used only in the companion paper \cite{Moore2}. We give it for completeness.

\bd \label{MooreFlow}
A {\rm Moore flow} is a $\mathcal{G}$-flow where $\mathcal{G}$ is the reparametrization category of Proposition~\ref{paramG}.
\ed

We now introduce some notations and definitions to provide more details. A $\mathcal{P}$-flow $X$ consists of a \textit{set of states} $X^0$, for each pair $(\alpha,\beta)$ of states a $\mathcal{P}$-space $\P_{\alpha,\beta}X$ of $\topdgr_0$ and for each triple $(\alpha,\beta,\gamma)$ of states an associative composition law \[*:\P_{\alpha,\beta}X \ot \P_{\beta,\gamma}X \to \P_{\alpha,\gamma}X,\]
i.e. for every triple of states $(\alpha,\beta,\gamma)\in X^0 \p X^0 \p X^0$, there is the commutative diagram of $\mathcal{P}$-spaces 
\[
\xymatrix@C=5em@R=5em
{
\P_{\alpha,\beta} X \ot \P_{\beta,\gamma} X \ot \P_{\gamma,\delta} X\fd{(\id,*)} \fr{(*,\id)} & \P_{\alpha,\gamma} X \ot \P_{\gamma,\delta} X \ar@{->}[d]^-{*}\\
\P_{\alpha,\beta} X \ot \P_{\beta,\delta} X \fr{*} & \P_{\alpha,\delta} X.
}
\]
A map of $\mathcal{P}$-flows $f$ from $X$ to $Y$ consists of a set map \[f^0:X^0 \to Y^0\] (often denoted by $f$ as well if there is no possible confusion) together for each pair of states $(\alpha,\beta)$ of $X$ with a natural transformation \[\P f:\P_{\alpha,\beta}X \longrightarrow \P_{f(\alpha),f(\beta)}Y\] such that the following diagram of $\topdgr_0$
\[
\xymatrix@C=5em@R=5em
{
	\P_{\alpha,\beta}X\ot \P_{\beta,\gamma}X \ar@{->}[r]^-{*}\ar@{->}[d] & \P_{\alpha,\gamma}X\ar@{->}[d]\\
	\P_{f(\alpha),f(\beta)}Y\ot \P_{f(\beta),f(\gamma)}Y \ar@{->}[r]^-{*} & \P_{f(\alpha),f(\gamma)}Y
}
\]
is commutative for all triples of states $(\alpha,\beta,\gamma)$ of $X$. 

The topological space $\P_{\alpha,\beta}X(\ell)$ is denoted by $\P_{\alpha,\beta}^\ell X$ and is called the space of \textit{execution paths of length $\ell$}. There is the sequence of homeomorphisms 
\[\begin{aligned}
\topdgr(\P_{\alpha,\beta}X &\ot \P_{\beta,\gamma}X, \P_{\alpha,\gamma}X)\\ & \iso \int_{(\ell_1,\ell_2)}\topdgr(\mathcal{P}(-,\ell_1+\ell_2)\p \P_{\alpha,\beta}^{\ell_1} X \p \P_{\beta,\gamma}^{\ell_2} X,\P_{\alpha,\gamma}X)\\
&\iso \int_{(\ell_1,\ell_2)}\ttop(\P_{\alpha,\beta}^{\ell_1} X \p \P_{\beta,\gamma}^{\ell_2} X,\P_{\alpha,\gamma}^{\ell_1+\ell_2}X)\\
&\iso [\mathcal{P}^{op}\p \mathcal{P}^{op},\top]\big((\ell_1,\ell_2)\mapsto \P_{\alpha,\beta}^{\ell_1} X \p \P_{\beta,\gamma}^{\ell_2} X, (\ell_1,\ell_2)\mapsto \P_{\alpha,\gamma}^{\ell_1+\ell_2}X\big),
\end{aligned}\]
the first homeomorphism by definition of $\ot$ and since $\topdgr$ is enriched, the second homeomorphism by \eqref{ey} and the third homeomorphism by \eqref{en}. 

Consequently, a $\mathcal{P}$-flow consists of a set of states $X^0$, for each pair $(\alpha,\beta)$ of states a $\mathcal{P}$-space $\P_{\alpha,\beta}X$ of $\topdgr_0$ and for each triple $(\alpha,\beta,\gamma)$ of states an associative composition law $*:\P^{\ell_1}_{\alpha,\beta}X \p \P^{\ell_2}_{\beta,\gamma}X \to \P^{\ell_1+\ell_2}_{\alpha,\gamma}X$ which is natural with respect to $(\ell_1,\ell_2)$. In other terms, for any map $\phi_i:\ell'_i\to \ell_i$ of $\mathcal{P}$ with $i=1,2$, there is the commutative diagram of topological spaces 
\[
\xymatrix@C=4em@R=4em
{
\P^{\ell_1}_{\alpha,\beta}X \p \P^{\ell_2}_{\beta,\gamma}X \fr{*}\ar@{->}[d]_-{(\P_{\alpha,\beta}(\phi_1),\P_{\beta,\gamma}(\phi_2))} & \P^{\ell_1+\ell_2}_{\alpha,\gamma}X \ar@{->}[d]^-{\P_{\alpha,\gamma}(\phi_1\ot \phi_2)}\\
\P^{\ell'_1}_{\alpha,\beta}X \p \P^{\ell'_2}_{\beta,\gamma}X \fr{*} & \P^{\ell'_1+\ell'_2}_{\alpha,\gamma}X
}
\] 
It means that for all $(x,y)\in \P^{\ell_1}_{\alpha,\beta}X \p \P^{\ell_2}_{\beta,\gamma}X$, there is the equality 
\[
(x.\phi_1) * (y.\phi_2) = (x*y).(\phi_1\ot \phi_2).
\]
The associativity means that 
\[\forall(\ell_1,\ell_2,\ell_3)\forall(\alpha,\beta,\gamma,\delta)\forall(x,y,z)\in  \P_{\alpha,\beta}^{\ell_1} X \p \P_{\beta,\gamma}^{\ell_2} X \p \P_{\gamma,\delta}^{\ell_3} X,(x*y)*z=x*(y*z).\]

\begin{prop}
	Every set $S$ gives rise to a $\mathcal{P}$-flow denoted by $S^\flat$ with the set of states $S$ and with the $\mathcal{P}$-space  $\Delta_{\mathcal{P}_0}\varnothing$ for each pair of states and to a $\mathcal{P}$-flow denoted by $S^\sharp$ with the set of states $S$ and with the $\mathcal{P}$-space $\Delta_{\mathcal{P}_0}\{0\}$ for each pair of states.
\end{prop}

\bpf The constant functors $\Delta_{\mathcal{P}_0}\varnothing$ and $\Delta_{\mathcal{P}_0}\{0\}$ belong to $\topdgr_0$ since $\id_\varnothing$ and $\id_{\{0\}}$ are continuous.
\epf

\bp \label{sk} The forgetful functor $(-)^0:\dtopP \to \set$ which takes a $\mathcal{P}$-flow $X$ to its set of state $X^0$ is both a left adjoint and a right adjoint.
\ep

\bpf
Let $X$ be a $\mathcal{P}$-flow. Let $S$ be a set. Then there are the bijections of sets $\set(X^0,S) \iso \dtopP(X,S^\sharp)$ and $\set(S,X^0) \iso \dtopP(S^\flat,X)$. 
\epf

\begin{defn}
	Let $X$ be a $\mathcal{P}$-flow. The {\rm $\mathcal{P}$-space of execution paths} $\P X$ of $X$ is by definition the $\mathcal{P}$-space \[\P X := \bigsqcup_{(\alpha,\beta)\in X^0\p X^0} \P_{\alpha,\beta}X.\]
	It yields a well-defined functor $\P:\dtopP \to \topdgr_0$. The image of $\ell$ is denoted by $\P^\ell U$. We therefore have the equality \[\P^\ell X = \bigsqcup_{(\alpha,\beta)\in X^0\p X^0} \P_{\alpha,\beta}^{\ell}X.\]
\end{defn}

\begin{nota}
	Let $X$ be a $\mathcal{P}$-flow. Denote by $\#X$ the {\rm cardinal of $X$}, i.e. the sum of the cardinal of the set of states of $X$ and of the cardinals of all topological spaces $\P_{\alpha,\beta}^{\ell}X$ for $\alpha,\beta$ running over $X^0$ and for $\ell$ running over the set of objects of $\mathcal{P}$.
\end{nota}

We want to emphasize the following elementary fact: 

\begin{lem} \label{emph-fact}
	Let $\C$ be a complete category. Let $X:I\to \C$ and $Y:I\to \C$ be two small diagrams of $\C$. Then the natural map $\limproj (X\p Y) \to \limproj X \p \limproj Y$ is an isomorphism.
\end{lem}

\bpf The functor $\limproj:\C^I \to \C$ is a right adjoint, the left adjoint being the constant diagram functor. Therefore, it preserves limits, and in particular binary product.
\epf

\bth \label{Mdtop-bicomplete}
The category $\dtopP$ is bicomplete.
\eth

\bpf
Let $X:I\longrightarrow \dtopP$ be a functor from a small category $I$ to $\dtopP$. Let $Y$ be the $\mathcal{P}$-flow defined as follows:
\begin{itemize}[leftmargin=*]
	\item The set of states $Y^0$ of $Y$ is defined as being the limit as sets $\limproj X(i)^0$ (Proposition~\ref{sk}).
	\item Let $\alpha,\beta\in \limproj X(i)^0$ and let $\alpha_i$ ($\beta_i$ resp.) be the image of $\alpha$ ($\beta$ resp.) in $X(i)^0$. Then let $\P_{\alpha,\beta} Y:= \limproj \P_{\alpha_i,\beta_i} X(i)$ where the limit is taken in the functor category $\topdgr_0$. There is therefore the isomorphism $\P_{\alpha,\beta}^{\ell} Y\iso \limproj \P_{\alpha_i,\beta_i}^{\ell} X(i)$ for all $\ell\in \Obj(\mathcal{P})$ and for all pairs $(\alpha,\beta)\in X^0\p X^0$.
	\item For $\alpha,\beta,\gamma\in \limproj X(i)^0$, let $\alpha_i$ ($\beta_i$, $\gamma_i$ resp.) be the image of $\alpha$ ($\beta$, $\gamma$ resp.) in $X(i)^0$. Then the composition map $*:\P_{\alpha,\beta}^\ell Y\p  \P_{\beta,\gamma}^{\ell'}Y\to \P_{\alpha,\gamma}^{\ell+\ell'}Y$ is taken as the limit of the $*_i:\P_{\alpha_i,\beta_i}^{\ell}X(i)\p \P_{\beta_i,\gamma_i}^{\ell'}X(i)\to \P_{\alpha_i,\gamma_i}^{\ell+\ell'}X(i)$ (we implicitly use Lemma~\ref{emph-fact}).
\end{itemize}
We obtain a well-defined $\mathcal{P}$-flow $Y$. It is clearly the limit $\limproj X$ in $\dtopP$. The constant diagram functor $\Delta_I:\dtopP\to \dtopP^I$ commutes with limits since limits in $\dtopP^I$ are calculated objectwise. Every map of $\dtopP^I$ $f:X\to \Delta_I Y$ certainly factors as a composite \[f:X\longrightarrow \Delta_I Z\longrightarrow \Delta_I Y\] with $\#Z < \sup_{i\in I} \#X_i$. There is a set of such $Z$ up to isomorphisms. We therefore have obtained a set of solutions. By Freyd's Adjoint Functor Theorem, the constant diagram functor $\Delta_I$ from $\dtopP$ to the category $\dtopP^I$ has a left adjoint. 
\epf

\begin{nota} \label{not-glob}
	Let $D$ be a $\mathcal{P}$-space. We denote by $\globP(D)$ the $\mathcal{P}$-flow defined as follows: 
	\[
	\begin{cases}
	\globP(D)^0 = \{0,1\}\\
	\P_{0,0}\globP(D)=\P_{1,1}\globP(D)=\P_{1,0}\globP(D)=\Delta_{\mathcal{P}_0}\varnothing\\
	\P_{0,1}\globP(D)=D
	\end{cases}	
	\]
	There is no composition law. This construction yields a functor \[\globP:\topdgr_0\to \dtopP.\]
\end{nota}

\bp \label{map-from-glob} Let $D$ be a $\mathcal{P}$-space. Let $X$ be a $\mathcal{P}$-flow. Then there is the natural bijection \[\dtopP(\globP(D),X)\iso \bigsqcup_{(\alpha,\beta)\in X^0\p X^0} \topdgr_0(D,\P_{\alpha,\beta}X).\]
\ep

\bpf
A map of $\mathcal{P}$-flows from $\globP(D)$ to $X$ is determined by the choice of two states $\alpha$ and $\beta$ of $X$ and by a map from $D$ to $\P_{\alpha,\beta}X$ because there is no composition law in $\globP(D)$.
\epf

\bth \label{dtopM-locpre} The category $\dtopP$ is locally presentable.
\eth

The category of small categories enriched over a closed monoidal category $\mathcal{V}$ is locally presentable as soon as $\mathcal{V}$ is locally presentable.  It is \cite[Theorem~4.5]{VCatLocPre} whose proof can probably be adapted to our situation: $\mathcal{P}$-flows are small semicategories enriched over a biclosed semimonoidal category which is locally presentable. It is shorter to proceed as follows.

\bpf
The category $\topdgr_0$ is locally presentable. Consider a dense generator $(G_i)_{i\in I}$ of $\lambda$-presentable $\mathcal{P}$-spaces for some regular cardinal $\lambda$. Let $f:X\to Y$ be a map of $\mathcal{P}$-flows. If $f$ induces a bijection from $\dtopP(\{0\}^\flat,X)=X^0$ to $\dtopP(\{0\}^\flat,Y)=Y^0$, then $f$ induces a bijection between the sets of states of $X$ and $Y$. If $f$ induces a bijection from $\dtopP(\globP(G_i),X)$ to $\dtopP(\globP(G_i),Y)$ for all $i\in I$, then by Proposition~\ref{map-from-glob}, $f$ induces an isomorphism between the $\mathcal{P}$-spaces $\P_{\alpha,\beta}X$ and $\P_{f(\alpha),f(\beta)}Y$ for all states $\alpha$ and $\beta$ of $X$ since the family $(G_i)_{i\in I}$ is a dense generator of $\topdgr_0$ by hypothesis. Thus the set of $\mathcal{P}$-flows \[\{\{0\}^\flat\} \cup \{\globP(G_i)\mid i\in I\}\] is a strong generator of $\dtopP$ by  \cite[Corollary~4.5.11]{Borceux1}. All $\mathcal{P}$-flows of \[\{\{0\}^\flat\} \cup \{\globP(G_i)\mid i\in I\}\] are $\lambda$-presentable. Therefore by \cite[Theorem~1.20]{TheBook} and by Theorem~\ref{Mdtop-bicomplete}, the category $\dtopP$ is locally $\lambda$-presentable.
\epf

\bp \label{adj-cp}
Let $D$ be a path-connected $\mathcal{P}$-space with $D\neq \Delta_{\mathcal{P}}\varnothing$. Then for every $\mathcal{P}$-flow $U$, there is a natural bijection 
\[
\topdgr_0(D,\P U) \iso \dtopP(\globP(D),U).
\]
\ep

\bpf
Let $f\in \topdgr_0(D,\P U)$. We have \[\P U = \bigsqcup_{(\alpha,\beta)\in U^0\p U^0} \P_{\alpha,\beta} U.\] For every $\ell\in \Obj(\mathcal{P})$, there is a continuous map $f_\ell:D(\ell) \to \P^\ell U$ which is natural with respect to $\ell$. For every $\ell\in \Obj(\mathcal{P})$, since $D(\ell)$ is path-connected and nonempty, the continuous map $f_\ell$ factors uniquely as a composite 
\[f_\ell:D(\ell) \longrightarrow \P_{\alpha_\ell,\beta_\ell}^\ell U \longrightarrow  \P^\ell U.\]
Consider a map $\phi:\ell'\to \ell$ of $\mathcal{P}$. It yields the commutative diagram of spaces 
\[
\xymatrix@C=3em@R=3em
{ &  \P_{\alpha_\ell,\beta_\ell}^\ell U \ar@{->}[rd]^-{\subset} & \\
	D(\ell) \ar@{->}[ru]\ar@{->}[rr]^-{f_\ell}\fd{D(\phi)}	 && \P^\ell U \ar@{->}[d]^-{\P U(\phi)} \\
	D(\ell') \ar@{->}[rd]\ar@{->}[rr]^-{f_{\ell'}}&& \P^{\ell'} U\\
	& \P_{\alpha_{\ell'},\beta_{\ell'}}^{\ell'} U \ar@{->}[ru]_-{\subset} &
}
\]
Since $D(\ell)$ is nonempty, there exists $x\in D(\ell)$. Then 
\[
\big(\P U)(\phi)\big)(f_\ell(x)) = f_{\ell'}(D(\phi)(x)) \in \P_{\alpha_\ell,\beta_\ell}^{\ell'}U \cap \P_{\alpha_{\ell'},\beta_{\ell'}}^{\ell'}U.
\] We deduce that \[(\alpha_\ell,\beta_\ell) = (\alpha_{\ell'},\beta_{\ell'}).\] Since the space $\mathcal{P}(\ell',\ell)$ is contractible for any $\ell,\ell'\in \Obj(\mathcal{P})$, it is nonempty. The mapping $\ell\mapsto (\alpha_\ell,\beta_\ell)$ from $\Obj(\mathcal{P})$ to $U^0\p U^0$ is therefore a constant which only depends on $f$. Consequently, the map of $\mathcal{P}$-spaces $f:D \to \P U$ factors uniquely as a composite 
\[f:D \longrightarrow \P_{\alpha,\beta} U \longrightarrow  \P U\]
for some $(\alpha,\beta)\in U^0 \p U^0$. This characterizes a map of $\mathcal{P}$-flow $\widehat{f}:\globP(D)\to U$. We easily see that the mapping $f\mapsto \widehat{f}$ is bijective. 
\epf

\bth \label{global-path-accessible} The path $\mathcal{P}$-space functor $\P:\dtopP\to \topdgr_0$ is  a right adjoint. In particular, it is limit-preserving and accessible. \eth

\bpf We mimick the construction made in \cite[Theorem~5.9]{leftproperflow} of the left adjoint of the path functor $\P:\dtop \to \top$. Using Proposition~\ref{cc-dgrtop}, write for a $\mathcal{P}$-space $D$: 
\[
D = \bigsqcup_{D_c\in \CC(D)} D_c.
\]
The left adjoint $\mathbf{G}:\topdgr_0 \to {\dtopP}$ is defined on objects by the formula
\[
\mathbf{G}(D) := \bigsqcup_{D_c\in \CC(D)} \globP(D_c).
\]
The definition of $\mathbf{G}:\topdgr_0 \to {\dtopP}$ on maps is clear. Choosing a map of flows from $\mathbf{G}(D)$ to a $\mathcal{P}$-flow $U$ is equivalent to choosing a map of $\mathcal{P}$-spaces from $D$ to $\P U$ because the image of any element of $\mathbf{G}(D)^0$ is forced by Proposition~\ref{adj-cp}. 
\epf

The composite functor 
\[
\xymatrix@C=3em
{
\dtopP \fr{\P} &  \topdgr_0 \fr{\subset} & \set^{\mathcal{P}_0}
}
\]
is finitely accessible since colimits are generated by free finite compositions of execution paths. On the contrary, it is unlikely that the functor $\P:\dtopP\to \topdgr_0$ is finitely accessible. However, in the case of a transfinite tower of objectwise relative-$T_1$ inclusions of spaces, one has: 

\bth  \label{global-path-almost-accessible}
Let $\lambda$ be a limit ordinal. Let $X:\lambda \to \dtopP$ be a colimit-preserving functor of $\mathcal{P}$-flows such that for all $\mu<\lambda$, the map of $\mathcal{P}$-spaces $\P X_\mu \to \P X_{\mu+1}$ is an objectwise relative-$T_1$ inclusion. Then the canonical map \[\liminj (\P.X) \longrightarrow \P \liminj X\] is an isomorphism of $\mathcal{P}$-spaces. Moreover the topology of $\P^{\ell} \liminj X$ is the final topology for all $\ell>0$.
\eth

\bpf
The topology of $\P^{\ell} \liminj X$ is the final topology for all $\ell>0$ because the relative-$T_1$ inclusions are one-to-one and because colimits are calculated objectwise in $\topdgr_0$. The rest of the proof is similar to the proof of \cite[Theorem~5.5]{leftproperflow}. The key point is to prove that the set map 
\[
\P^{\ell_1} \liminj X \p_{\liminj X^0}  \P^{\ell_2} \liminj X \longrightarrow \P^{\ell_1+\ell_2} \liminj X
\]
is continuous. It suffices to prove that every composite set map 
\[
[0,1] \to \P^{\ell_1} \liminj X \p_{\liminj X^0}  \P^{\ell_2} \liminj X \longrightarrow \P^{\ell_1+\ell_2} \liminj X
\]
is continuous as soon as the left-hand map is continuous. Since $[0,1]$ is finite relative to relative-$T_1$-inclusions by \cite[Proposition~2.5]{leftproperflow}, there exists an ordinal $\nu<\lambda$ such that the diagram
\[
\xymatrix@C=3em@R=3em
{
	[0,1] \fr{} \ar@{=}[d] & \P^{\ell_1} X_\nu \p_{X_\nu^0} \P^{\ell_2} X_\nu  \fr{*} \fd{} & \P^{\ell_1+\ell_2} X_\nu \fd{}\\
	[0,1] \fr{} & \P^{\ell_1} \liminj X \p_{\liminj X^0} \P^{\ell_2} \liminj X  \fr{*} & \P^{\ell_1+\ell_2} \liminj X
}
\]
is commutative. The top arrow $\P^{\ell_1} X_\nu \p_{X_\nu^0} \P^{\ell_2} X_\nu \to \P^{\ell_1+\ell_2} X_\nu$ is continuous because it is the composition law of a $\mathcal{P}$-flow. We deduce that the bottom arrow $\P^{\ell_1} \liminj X \p_{\liminj X^0} \P^{\ell_2} \liminj X  \to \P^{\ell_1+\ell_2} \liminj X$ is continuous as well by equipping $\P^{\ell} \liminj X$ with the final topology for all $\ell >0$.
\epf

\section{Reminder about model categories of fibrant objects}
\label{Isaev}

We start first by some notations: 

\begin{itemize}[leftmargin=*]
	\item $(-)^{cof}$ denotes a cofibrant replacement functor. 
	\item $f\boxslash g$ means that $f$ satisfies the \textit{left lifting property} (LLP) with respect to $g$, or equivalently when $g$ satisfies the \textit{right lifting property} (RLP) with respect to $f$. 
	\item $\inj(\C) = \{g \in \K, \forall f \in \C, f\boxslash g\}$.
	\item $\cof(\C)=\{f\mid \forall g\in \inj(\C), f\boxslash g\}$.
	\item $\cell(\C)$ is the class of transfinite compositions of pushouts of elements of $\C$.
\end{itemize}

All objects of all model categories of this paper are fibrant. We will be using the following characterization of a Quillen equivalence. A Quillen adjunction $F\dashv G:\C\leftrightarrows \D$ is a Quillen equivalence if and only if for all objects $X$ of $\D$, the natural map $F(G(X)^{cof})\to X$ is a weak equivalence of $\D$ (the functor is then said \textit{homotopically surjective}) and if for all cofibrant objects $Y$ of $\C$, the unit of the adjunction $Y\to G(F(Y))$ is a weak equivalence of $\C$.

We summarize in the following theorem what we want to use from \cite{Isaev}. For short, Isaev's paper gives a systematic way to construct model categories of fibrant objects. We already used some part of the following results in \cite{leftdetflow} to simplify the construction of the q-model category of flows. Unlike in \cite{leftdetflow}, we also want to add some comments about the class of weak equivalences of such model categories.

\bth \cite[Theorem~4.3, Proposition~4.4, Proposition~4.5 and Corollary~4.6]{Isaev} \label{Isa1}
Let $\K$ be a locally presentable category. Let $I$ be a set of maps of $\K$ such that the domains $D$ of the maps of $I$ are $I$-cofibrant (i.e. such that $\varnothing\to D$ belongs to $\cof(I)$). Suppose that for every map $i:U\to V \in I$, there exists an object $C_U(V)$ such that the relative codiagonal map $V\sqcup_U V \to V$ factors as a composite \[V\sqcup_U V \stackrel{\gamma_0\sqcup \gamma_1}\longrightarrow C_U(V)\longrightarrow V\] such that the left-hand map belongs to $\cof(I)$. Let \[J_I = \{\gamma_0:V\to C_U(V)\mid U\to V \in I\}.\] Suppose that there exists a path functor $\cocyl:\K \to \K$, i.e. an endofunctor of $\K$ equipped with two natural transformations $\tau:\id\Rightarrow\cocyl$ and $\pi:\cocyl\Rightarrow \id\p \id$ such that the composite $\pi.\tau$ is the diagonal. Moreover we suppose that the path functor satisfies the following hypotheses: 
\begin{enumerate}
	\item With $\pi=(\pi_0,\pi_1)$, $\pi_0:\cocyl(X)\to X$ and $\pi_1:\cocyl(X)\to X$ have the RLP with respect to $I$.
	\item The map $\pi:\cocyl(X)\to X\p X$ has the RLP with respect to the maps of $J_I$.
\end{enumerate}
Then there exists a unique model category structure on $\K$ such that the set of generating cofibrations is $I$ and such that the set of generating trivial cofibrations is $J_I$. Moreover, all objects are fibrant.
\eth

Isaev's paper contains also an interesting characterization of the class of weak equivalences which is recalled now: 

\bth \cite[Proposition~3.5]{Isaev} \label{Isa2}
With the notations of Theorem~\ref{Isa1}. A map $f:X\to Y$ of $\K$ is a weak equivalence if and only if it satisfies the RLP up to homotopy with respect to any map of $I$, i.e for any commutative diagram of solid arrows of $\K$ (where $i\in I$)
\[
\xymatrix@C=5em@R=5em
{
	U \ar@{->}[d]_-{i} \fr{u} & X \ar@{->}[d]^-{f}\\
	V \fr{v} \ar@{-->}[ru]^-{k}_-{\sim i} & Y
}
\]
there exists a lift $k$ such that $k.i=u$ and such that there exists a relative homotopy $h:C_U(V)\to Y$ between $f.k$ and $v$.
\eth

This characterization of the class of weak equivalences can be formulated in terms of injectivity class in the category of morphisms of $\K$ as follows: 

\bth \label{RLP-up-to-homotopy}
With the notations of Theorem~\ref{Isa1}. A map $f:X\to Y$ of $\K$ satisfies the RLP up to homotopy with respect to a map $U\to V$ of $I$ if and only if it is injective in $\Mor(\K)$ with respect to the map of maps 
\[
\arrow{U}{V} \longrightarrow \arrow{V}{C_U(V)}
\]
induced by the maps $V\sqcup_U V\to C_U(V)$ of $\K$.
\eth

\bpf
By Theorem~\ref{Isa2}, a map $f:X\to Y$ of $\K$ is a weak equivalence if and only if it satisfies the RLP up to homotopy with respect to any map of $I$. This definition can be reworded as follows. For any commutative diagram of solid arrows of $\K$ of the form 
\[
\xymatrix@C=5em@R=5em
{
	U \ar@{->}[r]^-{i} \ar@/^20pt/@{->}[rr]^-{u} \ar@{->}[d]_-{i} & V \fd{}\ar@{-->}[r]^-{k} & X \ar@{->}[d]^-{f} \\
	V \ar@/_20pt/@{->}[rr]_-{v} \fr{}& C_U(V)  \ar@{-->}[r]^-{h}  & Y,
}
\]
there exist $h$ and $k$ making the diagram commutative. It is exactly the injectivity condition when we regard this diagram as a diagram of $\Mor(\K)$. 
\epf

We can now reformulate Isaev's characterization of the class of weak equivalences of $\K$ as follows:

\begin{cor} \label{inj-charac-weak-equiv} 
	With the notations of Theorem~\ref{Isa1}. A map $f:X\to Y$ of $\K$ is a weak equivalence if and only if it is injective with respect to the maps of maps 
	\[
	\arrow{U}{V} \longrightarrow \arrow{V}{C_U(V)}
	\]
	induced by the map $V\sqcup_U V\to C_U(V)$ for $U\to V$ running over $I$. In particular, the class of weak equivalences of $\K$ is closed under small products because it is a small injectivity class of a locally presentable category.
\end{cor}

The idea of Corollary~\ref{inj-charac-weak-equiv} comes from a passing remark due to Jeff Smith and mentioned in \cite{weak-equiv-presheaves}. \cite[Theorem~14]{weakeq-algebraic} contains interesting results of the same kind for more general combinatorial model categories.

\section{Homotopy theory of \texorpdfstring{$\mathcal{P}$}{Lg}-flows}
\label{homotopytheoryMoore}

We are going to construct a model structure on $\dtopP$, called the \textit{q-model structure}, by mimicking the method used in \cite{leftdetflow} for the construction of the q-model structure of flows. 

We equip the category $\topdgr_0$ with the \textit{projective q-model structure} \cite[Theorem~6.2]{dgrtop}. It is denoted by $[\mathcal{P}^{op},\top_q]^{proj}_0$. The (trivial resp.) projective q-fibrations are the objectwise (trivial resp.) q-fibrations of spaces. The weak equivalences are the objectwise weak homotopy equivalences. The set of generating projective q-cofibrations is the set of maps \[\mathbb{I} = \{\mathbb{F}^{\mathcal{P}^{op}}_{\ell}\mathbf{S}^{n-1}\to \mathbb{F}^{\mathcal{P}^{op}}_{\ell}\mathbf{D}^n\mid n\geq 0,\ell\in \Obj(\mathcal{P})\}\] induced by the inclusions $\mathbf{S}^{n-1}\subset \mathbf{D}^n$. The set of generating trivial projective q-cofibra\-tions is the set of maps \[\mathbb{J} =\{ \mathbb{F}^{\mathcal{P}^{op}}_{\ell}\mathbf{D}^{n}\to \mathbb{F}^{\mathcal{P}^{op}}_{\ell}\mathbf{D}^{n+1}\mid n\geq 0,\ell\in \Obj(\mathcal{P})\}\] where the maps $\mathbf{D}^{n}\subset \mathbf{D}^{n+1}$ are induced by the mappings $(x_1,\dots,x_n) \mapsto (x_1,\dots,x_n,0)$.

\bp\label{ev-adj} \cite[Proposition~5.5 and Proposition~6.7]{dgrtop}
For every $\mathcal{P}$-space $F:\mathcal{P}^{op}\to \top$, every $\ell\in \Obj(\mathcal{P})$ and every topological space $X$, we have 
the natural bijection of sets \[\topdgr_0(\mathbb{F}^{\mathcal{P}^{op}}_{\ell}X,F) \iso \top(X,F(\ell)).\] In particular, the functor \[\mathbb{F}^{\mathcal{P}^{op}}_{\ell}:\top \to \topdgr_0\] is colimit-preserving for all $\ell\in \Obj(\mathcal{P})$. It induces a left Quillen functor $\top_q\to[\mathcal{P}^{op},\top_q]^{proj}_0$.
\ep

\begin{nota}  
	Let \[\begin{aligned}
	&I^{\mathcal{P}}=\{\globP(\mathbb{F}^{\mathcal{P}^{op}}_{\ell}\mathbf{S}^{n-1}) \subset \globP(\mathbb{F}^{\mathcal{P}^{op}}_{\ell}\mathbf{D}^{n})\mid n\geq 0,\ell \in \Obj(\mathcal{P})\}, \\
	&I^{\mathcal{P}}_+= I^{\mathcal{P}} \cup \{C:\varnothing^\flat \to \{0\}^\flat,R:\{0,1\}^\flat \to \{0\}^\flat\},\\
	&J^{\mathcal{P}} = \{\globP(\mathbb{F}^{\mathcal{P}^{op}}_{\ell}\mathbf{D}^{n}) \subset \globP(\mathbb{F}^{\mathcal{P}^{op}}_{\ell}\mathbf{D}^{n+1})\mid n\geq 0,\ell \in \Obj(\mathcal{P})\}.
	\end{aligned}\] 
\end{nota}

\bp \label{ortho3}  Let $U\to V$ be a continuous map. A morphism of $\mathcal{P}$-flows $f:X\to Y$ satisfies the RLP with respect to $\globP(\mathbb{F}^{\mathcal{P}^{op}}_{\ell}U)\to \globP(\mathbb{F}^{\mathcal{P}^{op}}_{\ell}V)$ for all $\ell\in \Obj(\mathcal{P})$ if and only if for any $\alpha,\beta\in X^0$ and any $\ell\in \Obj(\mathcal{P})$, the map of topological spaces $\P_{\alpha,\beta}^{\ell} X\to \P_{f(\alpha),f(\beta)}^{\ell}Y$ satisfies the RLP with respect to the continuous map $U\to V$. \ep

\bpf 
A morphism of $\mathcal{P}$-flows $f:X\to Y$ satisfies the RLP with respect to \[\globP(\mathbb{F}^{\mathcal{P}^{op}}_{\ell}U)\to \globP(\mathbb{F}^{\mathcal{P}^{op}}_{\ell}V)\] for all $\ell\in \Obj(\mathcal{P})$ if and only if the set maps 
\begin{multline*}
\dtopP(\globP(\mathbb{F}^{\mathcal{P}^{op}}_{\ell}V),X) \longrightarrow \\
\dtopP(\globP(\mathbb{F}^{\mathcal{P}^{op}}_{\ell}U),X) \p_{\dtopP(\globP(\mathbb{F}^{\mathcal{P}^{op}}_{\ell}U),Y)} \dtopP(\globP(\mathbb{F}^{\mathcal{P}^{op}}_{\ell}V),Y)
\end{multline*}
are onto for all $\ell\in \Obj(\mathcal{P})$. By Proposition~\ref{map-from-glob}, it is equivalent to saying that the set maps 
\begin{multline*}
\topdgr_0(\mathbb{F}^{\mathcal{P}^{op}}_{\ell}V,\P_{\alpha,\beta}X) \longrightarrow \\
\topdgr_0(\mathbb{F}^{\mathcal{P}^{op}}_{\ell}U,\P_{\alpha,\beta}X) \p_{\topdgr_0(\mathbb{F}^{\mathcal{P}^{op}}_{\ell}U,\P_{f(\alpha),f(\beta)}Y)} \topdgr_0(\mathbb{F}^{\mathcal{P}^{op}}_{\ell}V,\P_{f(\alpha),f(\beta)}Y)
\end{multline*}
are onto for all $(\alpha,\beta)\in X^0\p X^0$ and all $\ell\in \Obj(\mathcal{P})$. By Proposition~\ref{ev-adj}, it is equivalent to saying that the set maps 
\[\top(V,\P_{\alpha,\beta}^{\ell}X) \longrightarrow 
\top(U,\P_{\alpha,\beta}^{\ell}X) \p_{\top(U,\P_{f(\alpha),f(\beta)}^{\ell}Y)} \top(V,\P_{f(\alpha),f(\beta)}^{\ell}Y)\]
are onto for all $(\alpha,\beta)\in X^0\p X^0$ and all $\ell\in \Obj(\mathcal{P})$. It is equivalent to saying that for any $\alpha,\beta\in X^0$ and any $\ell\in \Obj(\mathcal{P})$, the map of topological spaces $\P_{\alpha,\beta}^{\ell} X\to \P_{f(\alpha),f(\beta)}^{\ell}Y$ satisfies the RLP with respect to the continuous map $U\to V$.
\epf

\bp\label{ortho2} Let $f$ be a morphism of $\mathcal{P}$-flows. Then the following assertions are equivalent: 1) $f$ is bijective on states, 2) $f$ satisfies the RLP with respect to $R:\{0,1\}^\flat\longrightarrow \{0\}^\flat$ and $C:\varnothing^\flat\subset \{0\}^\flat$.
\ep

\bpf
The map $f$ is one-to-one on states if and only if it satisfies the RLP with respect to $R:\{0,1\}^\flat\longrightarrow \{0\}^\flat$ and is onto on states if and only if it satisfies the RLP with respect to $C:\varnothing^\flat\subset \{0\}^\flat$.
\epf

\bp \label{connected-colim-globM}
The globe functor $\globP:\topdgr_0\to \dtopP$ preserves connected colimits.
\ep

Note that the connectedness hypothesis is necessary. Indeed, $D$ and $E$ being two $\mathcal{P}$-spaces, the $\mathcal{P}$-flow $\globP(D\sqcup E)$ has two states whereas the $\mathcal{P}$-flow $\glob(D) \sqcup \glob(E)$ has four states.

\bpf Let $D:I\to \topdgr_0$ be a functor where $I$ is a connected small category. We obtain the sequence of natural bijections:
\[\begin{aligned}
\dtopP(&\glob(\liminj D_i),X) \\& \iso \bigsqcup_{(\alpha,\beta)\in X^0\p X^0} \topdgr_0(\liminj D_i,\P_{\alpha,\beta}X) & \hbox{by Proposition~\ref{map-from-glob}}\\
&\iso \bigsqcup_{(\alpha,\beta)\in X^0\p X^0} \limproj\topdgr_0(D_i,\P_{\alpha,\beta}X)&\hbox{by definition of a (co)limit}\\ 
&\iso \limproj \bigsqcup_{(\alpha,\beta)\in X^0\p X^0} \topdgr_0(D_i,\P_{\alpha,\beta}X) & \hbox{by connectedness of $I$}\\
&\iso \limproj \dtopP(\glob(D_i),X) & \hbox{by Proposition~\ref{map-from-glob}}\\
&\iso \dtopP(\liminj \glob(D_i),X)&\hbox{by definition of a (co)limit}.
\end{aligned}\]
The proof is complete using the Yoneda lemma.
\epf

There is the obvious proposition: 

\bp
Let $X$ be a $\mathcal{P}$-flow. Let $U$ be a topological space. The following data assemble into a $\mathcal{P}$-flow denoted by $\{U,X\}_S$: 
\begin{itemize}[leftmargin=*]
	\item The set of states is $X^0$.
	\item $\P_{\alpha,\beta}\{U,X\}_S(\ell) := \ttop(U,\P_{\alpha,\beta}^\ell{X})$ for all $\ell>0$ and for all $(\alpha,\beta)\in X^0\p X^0$.
	\item The composition maps are defined by the composite maps 
	\[\ttop(U,\P_{\alpha,\beta}^{\ell}X)\p \ttop(U,\P_{\beta,\gamma}^{\ell'}X) \iso \ttop(U,\P_{\alpha,\beta}^{\ell}X\p \P_{\beta,\gamma}^{\ell'}X) \to \ttop(U,\P_{\alpha,\gamma}^{\ell+\ell'}X)\]
	for all $(\alpha,\beta,\gamma)\in X^0 \p X^0 \p X^0$ and all $\ell,\ell'\in \Obj(\mathcal{P})$. 
\end{itemize}
\ep

Note that we stick to the notations of \cite[Notation~7.6]{model3} and \cite[Notation~3.8]{leftdetflow}. Indeed, like for the case of flows, $\{U,X\}_S$ could be denoted by $X^U$ only if $U$ is connected and nonempty. The correct definition of $X^U$ in the non-connected case is as follows: 

\bd \label{exp}
Let $U$ be a topological space.  Let 
\[
X^U = \prod_{V\in\CC(U)} \{V,X\}_S
\]
The mapping $(U,X)\mapsto X^U$ induces a well-defined functor from $\top^{op}\p \dtopP$ to $\dtopP$.
\ed

We could prove that the functor $X\mapsto X^U$ has a left adjoint $X\mapsto X\ot U$ and that the axioms of tensored and cotensored categories are satisfied. One of the ingredients of the proof is that every $\Delta$-generated space is homeomorphic to the disjoint union of its nonempty connected components.

\bth \label{dtopM-model-structure}
There exists a unique model structure on $\dtopP$ such that $I^{\mathcal{P}}_+$ is the set of generating cofibrations and such that all objects are fibrant. The set of generating trivial cofibrations is $J^{\mathcal{P}}$. It is called {\rm the q-model structure}. 
\eth

\bpf
We are going to check all hypotheses of Theorem~\ref{Isa1}. The category $\dtopP$ is locally presentable by Theorem~\ref{dtopM-locpre}. By Proposition~\ref{ev-adj}, the map of $\mathcal{P}$-spaces \[\varnothing = \mathbb{F}^{\mathcal{P}^{op}}_{\ell}\varnothing \to \mathbb{F}^{\mathcal{P}^{op}}_{\ell}\mathbf{S}^{n-1}\] is a projective q-cofibration for all $\ell\in \Obj(\mathcal{P})$ and all $n\geq 0$. Thus all domains of all maps of $I^{\mathcal{P}}_+$ are cofibrant with respect to $I^{\mathcal{P}}_+$. We can factor the relative codiagonal map $\mathbf{D}^{n} \sqcup_{\mathbf{S}^{n-1}}\mathbf{D}^{n} \to \mathbf{D}^{n}$ as a composite \[\mathbf{D}^{n} \sqcup_{\mathbf{S}^{n-1}}\mathbf{D}^{n} \subset \mathbf{D}^{n+1} \longrightarrow \mathbf{D}^{n}\] for all $n\geq 0$. We obtain for all $\ell\in \Obj(\mathcal{P})$ a composite map of $\mathcal{P}$-spaces
\[
\xymatrix
{
	\mathbb{F}^{\mathcal{P}^{op}}_{\ell}(\mathbf{D}^{n} \sqcup_{\mathbf{S}^{n-1}}\mathbf{D}^{n}) \fr{} & \mathbb{F}^{\mathcal{P}^{op}}_{\ell}\mathbf{D}^{n+1} \fr{} & \mathbb{F}^{\mathcal{P}^{op}}_{\ell}\mathbf{D}^{n}
}
\]
such that the left-hand map is a projective q-cofibration by Proposition~\ref{ev-adj}. Thus for $U\to V$ being one of the maps $\globP(\mathbb{F}^{\mathcal{P}^{op}}_{\ell}\mathbf{S}^{n-1}) \subset \globP(\mathbb{F}^{\mathcal{P}^{op}}_{\ell}\mathbf{D}^{n})$ for $n\geq 0$ and $\ell\in \Obj(\mathcal{P})$, we set \[C_U(V) = \globP(\mathbb{F}^{\mathcal{P}^{op}}_{\ell}\mathbf{D}^{n+1}).\] For all $n\geq 0$, we have a pushout diagram of topological spaces
\[\xymatrix
{
	\mathbf{S}^{n}\fr{\iso}\fd{} & \mathbf{D}^{n} \sqcup_{\mathbf{S}^{n-1}}\mathbf{D}^{n}\fd{}\\
	\mathbf{D}^{n+1}\fr{} & \cocartesien\mathbf{D}^{n+1}
}\]
which gives rise to the pushout diagram of $\mathcal{P}$-spaces
\[\xymatrix
{
	\mathbb{F}^{\mathcal{P}^{op}}_{\ell}\mathbf{S}^{n}\fr{}\fd{} & \mathbb{F}^{\mathcal{P}^{op}}_{\ell}\mathbf{D}^{n} \sqcup_{\mathbb{F}^{\mathcal{P}^{op}}_{\ell}\mathbf{S}^{n-1}}\mathbb{F}^{\mathcal{P}^{op}}_{\ell}\mathbf{D}^{n}\fd{}\\
	\mathbb{F}^{\mathcal{P}^{op}}_{\ell}\mathbf{D}^{n+1}\fr{} & \cocartesien\mathbb{F}^{\mathcal{P}^{op}}_{\ell}\mathbf{D}^{n+1}
}\]
for all $\ell\in \Obj(\mathcal{P})$ by Proposition~\ref{ev-adj}. The latter diagram gives rise to the pushout diagram of $\mathcal{P}$-flows 
\[\xymatrix
{
	\globP(\mathbb{F}^{\mathcal{P}^{op}}_{\ell}\mathbf{S}^{n})\fd{}\fr{} & \globP(\mathbb{F}^{\mathcal{P}^{op}}_{\ell}\mathbf{D}^{n}) \sqcup_{\globP(\mathbb{F}^{\mathcal{P}^{op}}_{\ell}\mathbf{S}^{n-1})}\globP(\mathbb{F}^{\mathcal{P}^{op}}_{\ell}\mathbf{D}^{n})\fd{}\\
	\globP(\mathbb{F}^{\mathcal{P}^{op}}_{\ell}\mathbf{D}^{n+1})\fr{} & \cocartesien \globP(\mathbb{F}^{\mathcal{P}^{op}}_{\ell}\mathbf{D}^{n+1})
}\] 
for all $n\geq 0$ and all $\ell\in \Obj(\mathcal{P})$ by Proposition~\ref{connected-colim-globM}. This implies that for $U\to V$ being one of the maps \[\globP(\mathbb{F}^{\mathcal{P}^{op}}_{\ell}\mathbf{S}^{n-1}) \subset \globP(\mathbb{F}^{\mathcal{P}^{op}}_{\ell}\mathbf{D}^{n})\] for $n\geq 0$ and $\ell\in \Obj(\mathcal{P})$, the map $V\sqcup_U V\longrightarrow C_U(V)$ belongs to $\cell(I^{\mathcal{P}}_+)$. The map $C:\varnothing^\flat\to\{0\}^\flat$ gives rise to the relative codiagonal map $\{0\}^\flat\sqcup \{0\}^\flat\to \{0\}^\flat$. Thus we set \[C_{\varnothing^\flat}(\{0\}^\flat)=\{0\}^\flat.\] In this case, the map $V\sqcup_U V\longrightarrow C_U(V)$ is $R:\{0,1\}^\flat\to\{0\}^\flat$ which belongs to $\cell(I^{\mathcal{P}}_+)$. The map $R:\{0,1\}^\flat\to \{0\}^\flat$ gives rise to the relative codiagonal map $\id_{\{0\}}$. Thus we set \[C_{\{0,1\}^\flat}(\{0\}^\flat)=\{0\}^\flat.\] In this case, the map $V\sqcup_U V\longrightarrow C_U(V)$ is $\id_{\{0\}^\flat}$ which belongs to $\cell(I^{\mathcal{P}}_+)$. The set of generating trivial cofibrations will be therefore the set of maps  \[\globP(\mathbb{F}^{\mathcal{P}^{op}}_{\ell}\mathbf{D}^{n}) \subset \globP(\mathbb{F}^{\mathcal{P}^{op}}_{\ell}\mathbf{D}^{n+1})\] for $n\geq 0$ and $\ell\in \Obj(\mathcal{P})$. The composite map $\{0,1\}\subset [0,1] \to \{0\}$ yields a natural composite map of $\mathcal{P}$-flows 
\[
\xymatrix@C=5em
{
	X \iso X^{\{0\}} \fr{\tau_X} & X^{[0,1]} \fr{(\pi_0,\pi_1)} & X^{\{0,1\}} \iso X\p X
}
\] 
which gives rise to the composite continuous map \[\P_{\alpha,\beta}^{\ell}X \to \ttop([0,1],\P_{\alpha,\beta}^{\ell}X) \to \P_{\alpha,\beta}^{\ell}X \p\P_{\alpha,\beta}^{\ell}X\] on the spaces of paths for all $(\alpha,\beta)\in X^0\p X^0$ and all $\ell\in \Obj(\mathcal{P})$. We have obtained a path object in the sense of Theorem~\ref{Isa1}. Since the maps $\pi_0$ and $\pi_1$ are bijective on states, they satisfy the RLP with respect to \[\{C:\varnothing^\flat \to \{0\}^\flat,R:\{0,1\}^\flat \to \{0\}^\flat\}\] by Proposition~\ref{ortho2}. By Proposition~\ref{ortho3}, the maps $\pi_0$ and $\pi_1$ satisfy the RLP with respect to the maps of $\mathcal{P}$-flows \[\globP(\mathbb{F}^{\mathcal{P}^{op}}_{\ell}\mathbf{S}^{n-1}) \subset \globP(\mathbb{F}^{\mathcal{P}^{op}}_{\ell}\mathbf{D}^{n})\] for all $n\geq 0$ and for all $\ell\in \Obj(\mathcal{P})$ if and only if the evaluation maps \[\ttop([0,1],\P_{\alpha,\beta}^{\ell}X) \rightrightarrows \P_{\alpha,\beta}^{\ell}X\]  on $0$ and $1$ satisfy the RLP with respect to the inclusion $\mathbf{S}^{n-1}\subset \mathbf{D}^{n}$ for all $n\geq 0$, for all $\ell\in \Obj(\mathcal{P})$ and for all $(\alpha,\beta)\in X^0\p X^0$, i.e. if and only if the evaluation maps \[\ttop([0,1],\P_{\alpha,\beta}^{\ell}X) \rightrightarrows \P_{\alpha,\beta}^{\ell}X\] are trivial q-fibrations for all $(\alpha,\beta)\in X^0\p X^0$ and all $\ell\in \Obj(\mathcal{P})$. The latter fact is a consequence of the fact that the q-model structure of $\top$ is cartesian monoidal and from the fact that the inclusions $\{0\}\subset [0,1]$ and $\{1\}\subset [0,1]$ are trivial q-cofibrations of $\top$. Finally, we have to check that the map $\pi:\cocyl(X)\to X\p X$ satisfies the RLP with respect to the maps \[\glob(\mathbb{F}^{\mathcal{P}^{op}}_{\ell}\mathbf{D}^{n}) \subset \glob(\mathbb{F}^{\mathcal{P}^{op}}_{\ell}\mathbf{D}^{n+1})\] for all $n\geq 0$ and for all $\ell\in \Obj(\mathcal{P})$. By Proposition~\ref{ortho3} again, it suffices to prove that the map \[\ttop([0,1],\P_{\alpha,\beta}^{\ell} X) \to \ttop(\{0,1\},\P_{\alpha,\beta}^{\ell} X) = \P_{\alpha,\beta}^{\ell}X\p \P_{\alpha,\beta}^{\ell}X\] is a q-fibration of topological spaces for all $(\alpha,\beta)\in X^0\p X^0$ and all $\ell\in \Obj(\mathcal{P})$. Since the q-model structure of $\top$ is cartesian monoidal, this comes from the fact that the inclusion $\{0,1\}\subset [0,1]$ is a q-cofibration of $\top$.
\epf

\bth \label{Moore-fibration} A map of $\mathcal{P}$-flows $f:X\to Y$ is a q-fibration if and only if for all $(\alpha,\beta)\in X^0\p X^0$, the map of $\mathcal{P}$-spaces $\P_{\alpha,\beta}X\to \P_{f(\alpha),f(\beta)}Y$ is a projective q-fibration of $[\mathcal{P}^{op},\top_q]^{proj}_0$. \eth

\bpf
A set of generating trivial cofibrations of the model category of $\mathcal{P}$-flows is given by the set $\{\globP(\mathbb{F}^{\mathcal{P}^{op}}_{\ell}\mathbf{D}^n)\to \globP(\mathbb{F}^{\mathcal{P}^{op}}_{\ell}\mathbf{D}^{n+1})\mid n\geq 0,\ell\in \Obj(\mathcal{P})\}$. By Proposition~\ref{ortho3}, a map of $\mathcal{P}$-flows $f:X\to Y$ is a q-fibration if and only if for any $\alpha,\beta\in X^0$ and any $\ell\in \Obj(\mathcal{P})$ the map of topological spaces $\P_{\alpha,\beta}^{\ell} X\to \P_{f(\alpha),f(\beta)}^{\ell}Y$ satisfies the RLP with respect to the continuous maps $\mathbf{D}^n\to \mathbf{D}^{n+1}$ for $n\geq 0$, i.e. if and only if the map of $\mathcal{P}$-spaces $\P_{\alpha,\beta}X\to \P_{f(\alpha),f(\beta)}Y$ is a projective q-fibration of $[\mathcal{P}^{op},\top_q]^{proj}_0$.
\epf

\bth \label{Moore-trivial-fibration} A map of $\mathcal{P}$-flows $f:X\to Y$ is a trivial q-fibration if and only if $f$ induces a bijection between the set of states of $X$ and $Y$ and for all $(\alpha,\beta)\in X^0\p X^0$, the map of $\mathcal{P}$-spaces $\P_{\alpha,\beta}X\to \P_{f(\alpha),f(\beta)}Y$ is a trivial projective q-fibration of $[\mathcal{P}^{op},\top_q]^{proj}_0$. \eth

\bpf 
A set of generating q-cofibrations of the model category of $\mathcal{P}$-flows is given by the set 
$
\{\globP(\mathbb{F}^{\mathcal{P}^{op}}_{\ell}\mathbf{S}^{n-1}) \subset \globP(\mathbb{F}^{\mathcal{P}^{op}}_{\ell}\mathbf{D}^{n})\mid n\geq 0,\ell \in \Obj(\mathcal{P})\} \cup \{C:\varnothing^\flat \to \{0\}^\flat,R:\{0,1\}^\flat \to \{0\}^\flat\}
$. By Proposition~\ref{ortho3} and Proposition~\ref{ortho2}, a map of $\mathcal{P}$-flows $f:X\to Y$ is a trivial q-fibration if and only if it induces a bijection between the set of states and for any $\alpha,\beta\in X^0$ and any $\ell\in \Obj(\mathcal{P})$ the map of topological spaces $\P_{\alpha,\beta}^{\ell} X\to \P_{f(\alpha),f(\beta)}^{\ell}Y$ satisfies the RLP with respect to the continuous maps $\mathbf{S}^{n-1}\to \mathbf{D}^{n}$ for $n\geq 0$, i.e. if and only if the map of $\mathcal{P}$-spaces $\P_{\alpha,\beta}X\to \P_{f(\alpha),f(\beta)}Y$ is a trivial projective q-fibration of $[\mathcal{P}^{op},\top_q]^{proj}_0$.
\epf

The following proposition is a generalization of Proposition~\ref{map-from-glob}.

\bp \label{map-from-glob3}
Let $I$ be a small category having an initial object $\varnothing$. For all objects $i\in I$, let $f_i:\varnothing\to i$ be the unique morphism of $I$. Let $P$ be an object of the functor category $(\topdgr_0)^I$ (i.e. $P$ is a diagram of $\mathcal{P}$-spaces). Let $X$ be an object of the functor category $\dtopP^I$. The functor $\globP:\topdgr_0\to \dtopP$ inducing a functor $\globP:(\topdgr_0)^I\to \dtopP^I$ denoted in the same way, there is the natural bijection of sets
\[
\dtopP^I(\globP(P),X) \iso \bigsqcup_{(\alpha,\beta)\in X_\varnothing^0\p X_\varnothing^0}(\topdgr_0)^I(P,\P (X)),
\]
where the notation $\P(X)$ means the diagram of $\mathcal{P}$-spaces defined on objects by the mapping 
\[
i \mapsto \P_{f_i(\alpha),f_i(\beta)}X_i
\]

\ep

\bpf All maps of the diagram of $\mathcal{P}$-flows $\globP(P)$ are bijective on states. Therefore, a map $\globP(P)\to X$ is determined by the image $(\alpha,\beta)\in X_\varnothing^0\p X_\varnothing^0$ of the states $0$ and $1$ of the $\mathcal{P}$-flow $\globP(P_\varnothing)$ and, once this choice is done, by a map $P\to \P.X$. 
\epf

\begin{cor} \label{map-from-glob2} Let $f:X\to Y$ be a map of $\mathcal{P}$-flows. Let $P\to Q$ be a map of $\topdgr_0$. Then there is the natural bijection of sets 
	{\small\[
		\Mor(\dtopP)\left(\arrow{\globP(P)}{\globP(Q)},\arrow{X}{Y}\right) \iso \bigsqcup_{(\alpha,\beta)\in X^0\p X^0}\Mor\left(\topdgr_0\right)\left(\arrow{P}{Q},\arrow{\P_{\alpha,\beta}X}{\P_{f(\alpha),f(\beta)}Y}\right)
		\]}
\end{cor}

\bpf It is a corollary of Proposition~\ref{map-from-glob3} with the $2$-object small category $\varnothing \to \mathbf{1}$. \epf

The following well-known proposition is an easy consequence of \cite[page 219 (2)]{MR1712872}. It is explicitly stated for helping purpose. 

\bp \label{adj-functor-category} Let $\mathcal{L}\dashv \mathcal{R}:\C\leftrightarrows\D$ be a categorical adjunction. Let $I$ be a small category. Then 
\begin{itemize}[leftmargin=*]
	\item The functor $\mathcal{L}:\C\to \D$ induces a functor still denoted by $\mathcal{L}$ from the functor category $\C^I$ to the functor category $\D^I$ which takes $X$ to the functor $\mathcal{L}.X$.
	\item The functor $\mathcal{R}:\D\to \C$ induces a functor still denoted by $\mathcal{R}$ from the functor category $\D^I$ to the functor category $\C^I$ which takes $Y$ to the functor $\mathcal{R}.Y$.
\end{itemize}
We obtain an adjunction $\mathcal{L}\dashv \mathcal{R}:\C^I\leftrightarrows\D^I$.
\ep

\bpf
One has the natural bijections
\[\begin{aligned}
\D^I(\mathcal{L}(X),Y) &\iso \int_{i\in I} \D(\mathcal{L}(X(i)),Y(i)) &\hbox{by \cite[page 219 (2)]{MR1712872}}\\
&\iso \int_{i\in I} \C(X(i),\mathcal{R}(Y(i))) & \hbox{by adjunction} \\ & \iso \C^I(X,\mathcal{R}(Y))&\hbox{by \cite[page 219 (2)]{MR1712872}.}
\end{aligned}\]
\epf

\bp \label{ortho3-up-to-homotopy}  Let $n\geq 0$. Let $\ell\in \Obj(\mathcal{P})$. A morphism of $\mathcal{P}$-flows $f:X\to Y$ satisfies the RLP up to homotopy with respect to $\globP(\mathbb{F}^{\mathcal{P}^{op}}_{\ell}\mathbf{S}^{n-1})\to \globP(\mathbb{F}^{\mathcal{P}^{op}}_{\ell}\mathbf{D}^n)$ if and only if for all $\alpha,\beta\in X^0$, the map of topological spaces $\P_{\alpha,\beta}^{\ell} X\to \P_{f(\alpha),f(\beta)}^{\ell}Y$ satisfies the RLP up to homotopy with respect to the continuous map $\mathbf{S}^{n-1}\subset \mathbf{D}^n$. \ep

\bpf 
For $U\to V$ being one of the maps of \[\bigg\{\globP(\mathbb{F}^{\mathcal{P}^{op}}_{\ell}\mathbf{S}^{n-1}) \subset \globP(\mathbb{F}^{\mathcal{P}^{op}}_{\ell}\mathbf{D}^{n})\mid n\geq 0,\ell\in \Obj(\mathcal{P})\bigg\},\] we have (see the proof of Theorem~\ref{dtopM-model-structure}) \[C_U(V) = \globP(\mathbb{F}^{\mathcal{P}^{op}}_{\ell}\mathbf{D}^{n+1}).\]
Therefore, by Theorem~\ref{RLP-up-to-homotopy}, a morphism of $\mathcal{P}$-flows $f:X\to Y$ satisfies the RLP up to homotopy with respect to the map  $\globP(\mathbb{F}^{\mathcal{P}^{op}}_{\ell}\mathbf{S}^{n-1})\to \globP(\mathbb{F}^{\mathcal{P}^{op}}_{\ell}\mathbf{D}^n)$ if and only if it injective with respect to the map of maps 
\[\arrow{\globP(\mathbb{F}^{\mathcal{P}^{op}}_{\ell}\mathbf{S}^{n-1})}{\globP(\mathbb{F}^{\mathcal{P}^{op}}_{\ell}\mathbf{D}^n)} \longrightarrow \arrow{\globP(\mathbb{F}^{\mathcal{P}^{op}}_{\ell}\mathbf{D}^n)}{\globP(\mathbb{F}^{\mathcal{P}^{op}}_{\ell}\mathbf{D}^{n+1})}\]
in the category $\Mor(\dtopP)$, in other terms if and only if the set map 
\[
\Mor(\dtopP)\left(\arrow{\globP(\mathbb{F}^{\mathcal{P}^{op}}_{\ell}\mathbf{D}^n)}{\globP(\mathbb{F}^{\mathcal{P}^{op}}_{\ell}\mathbf{D}^{n+1})},\arrow{X}{Y}\right) \to \Mor(\dtopP)\left(\arrow{\globP(\mathbb{F}^{\mathcal{P}^{op}}_{\ell}\mathbf{S}^{n-1})}{\globP(\mathbb{F}^{\mathcal{P}^{op}}_{\ell}\mathbf{D}^n)},\arrow{X}{Y}\right)
\]
induced by the above map of maps is onto. However, we have the sequences of bijections 
\[\begin{aligned}
&\Mor(\dtopP)\left(\arrow{\globP(\mathbb{F}^{\mathcal{P}^{op}}_{\ell}\mathbf{D}^n)}{\globP(\mathbb{F}^{\mathcal{P}^{op}}_{\ell}\mathbf{D}^{n+1})},\arrow{X}{Y}\right)\\
&\iso \bigsqcup\limits_{(\alpha,\beta)\in X^0\p X^0}\Mor\left(\topdgr_0\right)\left(\arrow{\mathbb{F}^{\mathcal{P}^{op}}_{\ell}\mathbf{D}^n}{{\mathbb{F}^{\mathcal{P}^{op}}_{\ell}\mathbf{D}^{n+1}}},\arrow{\P_{\alpha,\beta}X}{\P_{f(\alpha),f(\beta)}Y}\right)&\hbox{by Corollary~\ref{map-from-glob2}}\\
&\iso \bigsqcup\limits_{(\alpha,\beta)\in X^0\p X^0}\Mor\left(\top\right)\left(\arrow{\mathbf{D}^n}{{\mathbf{D}^{n+1}}},\arrow{\P_{\alpha,\beta}^{\ell}X}{\P_{f(\alpha),f(\beta)}^{\ell}Y}\right)&\hbox{by Proposition~\ref{adj-functor-category}}
\end{aligned}\]
and 
\[\begin{aligned}
&\Mor(\dtopP)\left(\arrow{\globP(\mathbb{F}^{\mathcal{P}^{op}}_{\ell}\mathbf{S}^{n-1})}{\globP(\mathbb{F}^{\mathcal{P}^{op}}_{\ell}\mathbf{D}^n)},\arrow{X}{Y}\right)\\
&\iso \bigsqcup\limits_{(\alpha,\beta)\in X^0\p X^0}\Mor\left(\topdgr_0\right)\left(\arrow{\mathbb{F}^{\mathcal{P}^{op}}_{\ell}\mathbf{S}^{n-1}}{{\mathbb{F}^{\mathcal{P}^{op}}_{\ell}\mathbf{D}^n}},\arrow{\P_{\alpha,\beta}X}{\P_{f(\alpha),f(\beta)}Y}\right)&\hbox{by Corollary~\ref{map-from-glob2}}\\
&\iso \bigsqcup\limits_{(\alpha,\beta)\in X^0\p X^0}\Mor\left(\top\right)\left(\arrow{\mathbf{S}^{n-1}}{{\mathbf{D}^n}},\arrow{\P_{\alpha,\beta}^{\ell}X}{\P_{f(\alpha),f(\beta)}^{\ell}Y}\right)&\hbox{by Proposition~\ref{adj-functor-category}.}
\end{aligned}\]
Therefore a morphism of $\mathcal{P}$-flows $f:X\to Y$ satisfies the RLP up to homotopy with respect to the map  $\globP(\mathbb{F}^{\mathcal{P}^{op}}_{\ell}\mathbf{S}^{n-1})\to \globP(\mathbb{F}^{\mathcal{P}^{op}}_{\ell}\mathbf{D}^n)$ if and only if the set map 
\[
\Mor\left(\top\right)\left(\arrow{\mathbf{D}^n}{{\mathbf{D}^{n+1}}},\arrow{\P_{\alpha,\beta}^{\ell}X}{\P_{f(\alpha),f(\beta)}^{\ell}Y}\right) \to \Mor\left(\top\right)\left(\arrow{\mathbf{S}^{n-1}}{{\mathbf{D}^n}},\arrow{\P_{\alpha,\beta}^{\ell}X}{\P_{f(\alpha),f(\beta)}^{\ell}Y}\right)
\]
induced by the map of $\Mor(\top)$ 
\[
\arrow{\mathbf{S}^{n-1}}{\mathbf{D}^n} \longrightarrow \arrow{\mathbf{D}^n}{\mathbf{D}^{n+1}}
\]
is onto. By Theorem~\ref{RLP-up-to-homotopy}, the latter condition is equivalent to saying that the map $\P_{\alpha,\beta}^{\ell}X\to \P_{f(\alpha),f(\beta)}^{\ell}Y$ satisfies the RLP up to homotopy with respect to the continuous map $\mathbf{S}^{n-1}\subset \mathbf{D}^n$.
\epf

\bp \label{up-to-homotopy-RC} A morphism of $\mathcal{P}$-flows $f:X\to Y$ satisfies the RLP up to homotopy with respect to $C:\varnothing^\flat \to \{0\}^\flat$ if and only if it is onto on states. A morphism of $\mathcal{P}$-flows $f:X\to Y$ satisfies the RLP up to homotopy with respect to $R:\{0,1\}^\flat\to\{0\}^\flat$ if and only if it is one-to-one on states. 
\ep

\bpf
Obvious.
\epf

\bth A map of $\mathcal{P}$-flows $f:X\to Y$ is a weak equivalence if and only if $f$ induces a bijection between the set of states of $X$ and $Y$ and for all $(\alpha,\beta)\in X^0\p X^0$, the map of $\mathcal{P}$-spaces $\P_{\alpha,\beta}X\to \P_{f(\alpha),f(\beta)}Y$ is a weak equivalence. \eth

\bpf By Theorem~\ref{Isa2}, a map of $\mathcal{P}$-flows $f:X\to Y$ is a weak equivalence if and only if it satisfies the RLP up to homotopy with respect to the maps of \begin{multline*}
I^{\mathcal{P}}_+=\{\globP(\mathbb{F}^{\mathcal{P}^{op}}_{\ell}\mathbf{S}^{n-1}) \subset \globP(\mathbb{F}^{\mathcal{P}^{op}}_{\ell}\mathbf{D}^{n})\mid n\geq 0,\ell \in \mathcal{P}\} \\\cup \{C:\varnothing^\flat \to \{0\}^\flat,R:\{0,1\}^\flat \to \{0\}^\flat\}.
\end{multline*} By Proposition~\ref{up-to-homotopy-RC} and by Proposition~\ref{ortho3-up-to-homotopy}, a map of $\mathcal{P}$-flows $f:X\to Y$ is then a weak equivalence if and only if $f$ induces a bijection between the set of states of $X$ and $Y$ and for all $(\alpha,\beta)\in X^0\p X^0$ and all $\ell \in \Obj(\mathcal{P})$, the map $\P_{\alpha,\beta}^{\ell}X\to \P_{f(\alpha),f(\beta)}^{\ell}Y$ satisfies the RLP up to homotopy with respect to the maps $\mathbf{S}^{n-1}\subset \mathbf{D}^n$ for all $n\geq 0$. By Theorem~\ref{Isa2} applied to the q-model category $\top$, it implies that a map of $\mathcal{P}$-flows $f:X\to Y$ is then a weak equivalence if and only if $f$ induces a bijection between the set of states of $X$ and $Y$ and for all $(\alpha,\beta)\in X^0\p X^0$ and all $\ell \in \Obj(\mathcal{P})$, the map $\P_{\alpha,\beta}^{\ell}X\to \P_{f(\alpha),f(\beta)}^{\ell}Y$ is a weak homotopy equivalence of topological spaces.
\epf

\section{Q-cofibrant \texorpdfstring{$\mathcal{P}$}{Lg}-flows have a projective q-cofibrant \texorpdfstring{$\mathcal{P}$}{Lg}-space of execution paths}
\label{cofibrant-moore-flow}

\bd
Let $f:D \to E$ and $g:F\to G$ be two maps of $\mathcal{P}$-spaces. Then the {\rm pushout product, with respect to $\ot$,} denoted by $f\boxtimes g$, of $f$ and $g$ is the map of $\mathcal{P}$-spaces
\[
f\boxtimes g:(D \ot G) \sqcup_{D\ot F} (E \ot F) \longrightarrow E\ot G.
\]
induced by the universal property of the pushout.
\ed

\begin{nota}
	$f\square g$ denotes the pushout product with respect to the binary product of two maps $f$ and $g$ of a bicomplete category. 
\end{nota}

\bp \label{precalcul0}
Let $f:U\to V$ and $f':U'\to V'$ be two maps of $\top$. Let $\ell,\ell'\in \Obj(\mathcal{P})$. Then there is the isomorphism of $\mathcal{P}$-spaces
\[
\mathbb{F}^{\mathcal{P}^{op}}_{\ell}(f) \boxtimes \mathbb{F}^{\mathcal{P}^{op}}_{\ell'}(g) \iso \mathbb{F}^{\mathcal{P}^{op}}_{\ell+\ell'}(f \square g).
\]
\ep 

\bpf
The domain of $\mathbb{F}^{\mathcal{P}^{op}}_{\ell}(f) \boxtimes \mathbb{F}^{\mathcal{P}^{op}}_{\ell'}(g)$ is the $\mathcal{P}$-space 
\[\begin{aligned}
\big(\mathbb{F}^{\mathcal{P}^{op}}_{\ell}(U) \ot \mathbb{F}^{\mathcal{P}^{op}}_{\ell'}(V')\big) \sqcup_{\mathbb{F}^{\mathcal{P}^{op}}_{\ell}(U) \ot \mathbb{F}^{\mathcal{P}^{op}}_{\ell}(U')} &\big(\mathbb{F}^{\mathcal{P}^{op}}_{\ell}(U') \ot \mathbb{F}^{\mathcal{P}^{op}}_{\ell'}(V)\big) \\ 
&\iso \mathbb{F}^{\mathcal{P}^{op}}_{\ell+\ell'}(U\p V') \sqcup_{\mathbb{F}^{\mathcal{P}^{op}}_{\ell+\ell'}(U\p V)} \mathbb{F}^{\mathcal{P}^{op}}_{\ell+\ell'}(U'\p V)\\
&\iso \mathbb{F}^{\mathcal{P}^{op}}_{\ell+\ell'} \big((U\p V') \sqcup_{(U\p V)}(U'\p V)\big),
\end{aligned}\]
the first isomorphism by Proposition~\ref{Ftenseur} and the second isomorphism by Proposition~\ref{ev-adj}. By Proposition~\ref{Ftenseur} again, the codomain of $\mathbb{F}^{\mathcal{P}^{op}}_{\ell}(f) \boxtimes \mathbb{F}^{\mathcal{P}^{op}}_{\ell'}(g)$ is 
\[
\mathbb{F}^{\mathcal{P}^{op}}_{\ell}(U') \ot \mathbb{F}^{\mathcal{P}^{op}}_{\ell'}(V') \iso \mathbb{F}^{\mathcal{P}^{op}}_{\ell+\ell'}(U'\p V'). 
\]
\epf

\bth  \label{modelsemimonoidal}
Let $f:D \to E$ and $g:F\to G$ be two maps of $\mathcal{P}$-spaces. Then 
\begin{itemize}[leftmargin=*]
	\item If $f$ and $g$ are projective q-cofibrations, then $f\boxtimes g$ is a projective q-cofibration.
	\item If moreover $f$ or $g$ is a trivial projective q-cofibration, then $f\boxtimes g$ is a trivial projective q-cofibration.
\end{itemize}
\eth

\bpf 
The projective q-cofibrations of $\mathcal{P}$-spaces are generated by the set of maps \[\mathbb{I} = \{\mathbb{F}^{\mathcal{P}^{op}}_{\ell}\mathbf{S}^{n-1}\to \mathbb{F}^{\mathcal{P}^{op}}_{\ell}\mathbf{D}^n\mid n\geq 0,\ell\in \Obj(\mathcal{P})\}\] induced by the inclusions $\mathbf{S}^{n-1}\subset \mathbf{D}^n$. The set of generating trivial projective q-cofibra\-tions is the set of maps \[\mathbb{J} = \{\mathbb{F}^{\mathcal{P}^{op}}_{\ell}\mathbf{D}^{n}\to \mathbb{F}^{\mathcal{P}^{op}}_{\ell}\mathbf{D}^{n+1}\mid n\geq 0,\ell\in \Obj(\mathcal{P})\}\] where the maps $\mathbf{D}^{n}\subset \mathbf{D}^{n+1}$ are induced by the mappings $(x_1,\dots,x_n) \mapsto (x_1,\dots,x_n,0)$. From Proposition~\ref{precalcul0}, Proposition~\ref{ev-adj} and from the fact that the q-model structure of $\top$ is cartesian monoidal, we deduce the inclusions 
\[\begin{aligned}
&\mathbb{I}\boxtimes \mathbb{I} \subset \cof(\mathbb{I}),\\
&\mathbb{I}\boxtimes \mathbb{J} \subset \cof(\mathbb{J}),\\
&\mathbb{J}\boxtimes \mathbb{I} \subset \cof(\mathbb{J}).
\end{aligned}\]
The proof is complete thanks to Theorem~\ref{closedsemi} and \cite[Lemma~4.2.4]{MR99h:55031}. 
\epf

Let $S$ be a nonempty set. Let $\mathcal{P}^{u,v}(S)$ be the small category defined by generators and relations as follows (see \cite[Section~3]{leftproperflow}): 
\begin{itemize}[leftmargin=*]
	\item $u,v\in S$ ($u$ and $v$ may be equal).
	\item The objects are the tuples of the form 
	\[\underline{m}=((u_0,\epsilon_1,u_1),(u_1,\epsilon_2,u_2),\dots ,(u_{n-1},\epsilon_n,u_n))\]
	with $n\geq 1$, $u_0,\dots,u_n \in S$ and \[\forall i\hbox{ such that } 1\leq i\leq n, \epsilon_i = 1\Rightarrow (u_{i-1},u_i)=(u,v).\] 
	\item There is an arrow \[c_{n+1}:(\underline{m},(x,0,y),(y,0,z),\underline{n}) \to (\underline{m},(x,0,z),\underline{n})\]
	for every tuple $\underline{m}=((u_0,\epsilon_1,u_1),(u_1,\epsilon_2,u_2),\dots ,(u_{n-1},\epsilon_n,u_n))$ with $n\geq 1$ and every tuple $\underline{n}=((u'_0,\epsilon'_1,u'_1),(u'_1,\epsilon'_2,u'_2),\dots ,(u'_{n'-1},\epsilon'_{n'},u'_{n'}))$ with $n'\geq 1$. It is called a \textit{composition map}. 
	\item There is an arrow \[I_{n+1}:(\underline{m},(u,0,v),\underline{n}) \to (\underline{m},(u,1,v),\underline{n})\] for every tuple $\underline{m}=((u_0,\epsilon_1,u_1),(u_1,\epsilon_2,u_2),\dots ,(u_{n-1},\epsilon_n,u_n))$ with $n\geq 1$ and every tuple $\underline{n}=((u'_0,\epsilon'_1,u'_1),(u'_1,\epsilon'_2,u'_2),\dots ,(u'_{n'-1},\epsilon'_{n'},u'_{n'}))$ with $n'\geq 1$.
	It is called an \textit{inclusion map}. 
	\item There are the relations (group A) $c_i.c_j = c_{j-1}.c_i$ if $i<j$ (which means since $c_i$ and $c_j$ may correspond to several maps that if $c_i$ and $c_j$ are composable, then there exist $c_{j-1}$ and $c_i$ composable satisfying the equality). 
	\item There are the relations (group B) $I_i.I_j = I_j.I_i$ if $i\neq j$. By definition of these maps, $I_i$ is never composable with itself. 
	\item There are the relations (group C) \[c_i.I_j = \begin{cases}
	I_{j-1}.c_i&\hbox{if } j\geq i+2\\
	I_j.c_i&\hbox{if } j\leq i-1.
	\end{cases}\]
	By definition of these maps, $c_i$ and $I_i$ are never composable as well as $c_i$ and $I_{i+1}$. 
\end{itemize}

By \cite[Proposition~3.7]{leftproperflow}, there exists a structure of Reedy category on $\mathcal{P}^{u,v}(S)$ with the $\mathbb{N}$-valued degree map defined by \[d((u_0,\epsilon_1,u_1),(u_1,\epsilon_2,u_2),\dots ,(u_{n-1},\epsilon_n,u_n)) = n + \sum_i \epsilon_i.\]
The maps raising the degree are the inclusion maps in the above sense. The maps decreasing the degree are the composition maps in the above sense.

\begin{nota}
	The latching object at $\underline{n}$ of a diagram $\D$ over $\mathcal{P}^{u,v}(S)$ is denoted by $L_{\underline{n}} \D$. 
\end{nota}

We recall the important theorem: 

\bth \cite[Theorem~3.9]{leftproperflow} \label{ok} Let $\K$ be a model category. Let $S$ be a nonempty set. Let $u,v\in S$. Let $\Dcat(\mathcal{P}^{u,v}(S),\K)$ be the category of functors and natural transformations from $\mathcal{P}^{u,v}(S)$ to $\K$. Then there exists a unique model structure on \[\Dcat(\mathcal{P}^{u,v}(S),\K)\] such that the weak equivalences are the objectwise weak equivalences and such that a map of diagrams $f:\D\to \mathcal{E}$ is a cofibration (called a Reedy cofibration) if for all objects $\underline{n}$ of $\mathcal{P}^{u,v}(S)$, the canonical map $L_{\underline{n}}\mathcal{E} \sqcup_{L_{\underline{n}} \D} \D(\underline{n})\to \mathcal{E}(\underline{n})$ is a cofibration of $\K$. Moreover the colimit functor $\liminj : \Dcat(\mathcal{P}^{u,v}(S),\K) \to \K$ is a left Quillen adjoint. \eth

Let $\de Z \to Z$ be a map of $\mathcal{P}$-spaces. Consider a pushout diagram of $\mathcal{P}$-flows 
\[
\xymatrix@C=4em@R=4em
{
	\glob(\partial Z) \fd{} \fr{g} & A \ar@{->}[d]^-{f} \\
	\glob(Z) \fr{\widehat{g}} & \cocartesien X.
}
\]
Let $T$ be the $\mathcal{P}$-space defined by the pushout diagram of $\topdgr_0$
\[
\xymatrix@C=4em@R=4em
{
	\partial Z  \fd{} \fr{g} & \P_{g(0),g(1)} A \ar@{->}[d]^-{f} \\
	Z  \fr{\widehat{g}} & \cocartesien T.
}
\]
Consider the diagram of $\mathcal{P}$-spaces $\D^f:\mathcal{P}^{g(0),g(1)}(A^0)\to \topdgr_0$ defined as follows:
\[
\D^f((u_0,\epsilon_1,u_1),(u_1,\epsilon_2,u_2),\dots ,(u_{n-1},\epsilon_n,u_n)) = Z_{u_0,u_1}\ot Z_{u_1,u_2} \ot \dots \ot Z_{u_{n-1},u_n}
\]
with 
\[
Z_{u_{i-1},u_i}=
\begin{cases}
\P_{u_{i-1},u_i}A & \hbox{if }\epsilon_i=0\\
T & \hbox{if }\epsilon_i=1
\end{cases}
\] 
In the case $\epsilon_i=1$, $(u_{i-1},u_i)=(g(0),g(1))$ by definition of $\mathcal{P}^{g(0),g(1)}(A^0)$. The inclusion maps $I_i's$ are induced by the map $f:\P_{g(0),g(1)} A \to T$. The composition maps $c_i's$ are induced by the composition law of $A$.

\bth \label{colim-calcul}
We obtain a well-defined diagram of $\mathcal{P}$-spaces \[\D^f:\mathcal{P}^{g(0),g(1)}(A^0)\to \topdgr_0.\] There is the isomorphism of $\mathcal{P}$-spaces $\liminj \D^f \iso \P X$. 
\eth

\bpf
A flow is a small semicategory enriched over the closed (semi)monoidal category $(\top,\p)$ (cf. Definition~\ref{def_flow}). A $\mathcal{P}$-flow is a small semicategory enriched over the biclosed semimonoidal category $(\topdgr_0,\ot)$. Therefore, the first assertion is proved like \cite[Proposition~4.6]{leftproperflow} and the second assertion like \cite[Theorem~4.8]{leftproperflow}.
\epf

We recall the explicit calculation of the pushout product of several morphisms, here in the bicomplete biclosed semimonoidal category $(\topdgr_0,\ot)$. 

\bp\label{calculpushout} 
Let $f_i:A_i\longrightarrow B_i$ for $0\leq i\leq p$ be $p+1$ maps of $\mathcal{P}$-spaces. Let $S\subset
\{0,\dots,p\}$. Let 
\[C_p(S):=C_0 \ot \dots \ot C_p \hbox{ with } \begin{cases}
C_i = A_i & \hbox{ if }i\notin S\\
C_i = B_i & \hbox{ if }i\in S.
\end{cases}\]
If $S$ and $T$ are two subsets of $\{0,\dots,p\}$ such that $S\subset
T$, let \[C_p(i_S^T):C_p(S)\longrightarrow C_p(T)\] be the morphism 
\[g_0 \ot \dots \ot g_p \hbox{ with } \begin{cases}
g_i = \id_{B_i} & \hbox{ if }i\in S\\
g_i = f_i & \hbox{ if }i\in T\backslash S\\
g_i = \id_{A_i} & \hbox{ if }i\notin T.
\end{cases}\]
Then:  
\begin{enumerate} 
	\item the mappings $S\mapsto C_p(S)$ and $i_S^T\mapsto C_p(i_S^T)$ 
	give rise to a functor from the order
	complex of the poset $\{0<\dots <p\}$ to $\topdgr_0$
	\item there exists a canonical morphism 
		\[\liminj_{S\subsetneqq \{0,\dots,p\}} 
		C_p(S)\longrightarrow C_p(\{0,\dots,p\}).\]
	and it is equal to the morphism $f_0\boxtimes\dots \boxtimes f_p$. 
\end{enumerate}
\ep 

\bpf
The proof is similar to the proof of \cite[Theorem~B.3]{3eme} for the case of $(\top,\p)$. It is still valid here because the semimonoidal structure $\ot$ is associative and biclosed by Theorem~\ref{closedsemi}. The statement has to be slightly modified since $\ot$ is not assumed to be symmetric.
\epf

\bp \label{prep1}
With the notations above.  Let $\underline{n}\in \Obj(\mathcal{P}^{g(0),g(1)}(A^0))$ with \[\underline{n} = ((u_0,\epsilon_1,u_1),(u_1,\epsilon_2,u_2),\dots ,(u_{n-1},\epsilon_n,u_n)).\] Then the continuous map \[L_{\underline{n}} \D^f \longrightarrow \D^f(\underline{n})\] is the pushout product $\boxtimes$ of the maps $\varnothing \to \P_{u_{i-1},u_i}A$ for $i$ running over $\{i\in [1,n]| \epsilon_i = 0\}$ and of the maps $\P_{g(0),g(1)} A \to T$ for $i$ running over $\{i\in [1,n]| \epsilon_i = 1\}$. Moreover, if for all $i\in [1,n]$, we have $\epsilon_i=0$, then $L_{\underline{n}} \D^f = \varnothing$. 
\ep

\bpf
It is a consequence of Proposition~\ref{calculpushout}. 
\epf

\bth \label{pretower}
With the notations above. Assume that the map of $\mathcal{P}$-spaces $\de Z \to Z$ is a (trivial resp.) projective q-cofibration. If $\P A$ is a projective q-cofibrant $\mathcal{P}$-space, then $\P f:\P A \to \P X$ is a (trivial resp.) projective q-cofibration of $\mathcal{P}$-spaces. 
\eth

\bpf
The particular case $\de Z = Z$, $f=\id_A$ and $A=X$ yields the isomorphism of $\mathcal{P}$-spaces \[\liminj \D^{\id_A} \iso \P A.\] We have a map of diagrams $\D^{\id_A} \to \D^{f}$ which induces for all $\underline{n}\in \Obj(\mathcal{P}^{g(0),g(1)}(A^0))$ a map of $\mathcal{P}$-spaces\[L_{\underline{n}} \D^f \sqcup_{L_{\underline{n}} \D^{\id_A}} \D^{\id_A}(\underline{n}) \longrightarrow \D^{f}(\underline{n}).\] Let \[\underline{n} = ((u_0,\epsilon_1,u_1),(u_1,\epsilon_2,u_2),\dots ,(u_{n-1},\epsilon_n,u_n)).\] There are two mutually exclusive cases: 
\begin{enumerate}
	\item[(a)] All $\epsilon_i$ for $i=1,\dots n$ are equal to zero.
	\item[(b)] There exists $i\in [1,n]$ such that $\epsilon_i = 1$.
\end{enumerate}
In the case (a), we have \[\D^{\id_A}(\underline{n}) = \D^{f}(\underline{n}) = \P_{u_0,u_1}A \ot \dots \ot \P_{u_{n-1},u_n}A.\] Moreover, by Proposition~\ref{prep1}, we have $L_{\underline{n}} \D^{\id_A} = L_{\underline{n}} \D^f = \varnothing$. We deduce that the map \[L_{\underline{n}} \D^f \sqcup_{L_{\underline{n}} \D^{\id_A}} \D^{\id_A}(\underline{n}) \longrightarrow \D^{f}(\underline{n})\] is isomorphic to the identity of $\D^{\id_A}(\underline{n})$. In the case (b), The map \[L_{\underline{n}} \D^{\id_A} \longrightarrow \D^{\id_A}(\underline{n})\] is by Proposition~\ref{prep1} a pushout product of several maps such that one of them is the identity map $\id: \P _{g(0),g(1)}A \to \P _{g(0),g(1)}A$ because $\epsilon_i=1$ for some $i$. Therefore the map $L_{\underline{n}} \D^{\id_A} \to \D^{\id_A}(\underline{n})$ is an isomorphism. We deduce that the map \[L_{\underline{n}} \D^f \sqcup_{L_{\underline{n}} \D^{\id_A}} \D^{\id_A}(\underline{n}) \longrightarrow \D^{f}(\underline{n})\] is isomorphic to the map $L_{\underline{n}} \D^f \to \D^f(\underline{n})$. By Proposition~\ref{prep1}, the map $L_{\underline{n}}\D^f \to \D^f(\underline{n})$ is a pushout product $\boxtimes$ of maps of the form $\varnothing\to \P _{\alpha,\beta}A$ and of the form $f:\P_{g(0),g(1)} A\to T$ for all objects $\underline{n}\in \Obj(\mathcal{P}^{g(0),g(1)}(A^0))$. We conclude that the map 
\[L_{\underline{n}} \D^f \sqcup_{L_{\underline{n}} \D^{\id_A}} \D^{\id_A}(\underline{n}) \longrightarrow \D^{f}(\underline{n})\]
is for all $\underline{n}$ either an isomorphism, or a pushout product of maps of the form $\varnothing\to \P _{\alpha,\beta}A$ and of the form $f:\P_{g(0),g(1)} A\to T$, the latter appearing at least once in the pushout product and being a pushout of the map of $\mathcal{P}$-spaces $\de Z\to Z$. We are now ready to complete the proof.

Suppose now that $\P A$ is a projective q-cofibrant $\mathcal{P}$-space. Therefore for all $(\alpha,\beta)\in A^0\p A^0$, the $\mathcal{P}$-space $\P_{\alpha,\beta} A$ is projective q-cofibrant. We deduce that the map $L_{\underline{n}}\D^f \to \D^f(\underline{n})$ is always a projective q-cofibration of $\topdgr_0$ for all $\underline{n}$ by Theorem~\ref{modelsemimonoidal}. We deduce that the map of diagrams $\D^{\id_A} \to \D^{f}$ is a Reedy projective q-cofibration. Therefore by passing to the colimit which is a left Quillen adjoint by Theorem~\ref{ok}, we deduce that the map $\P A \to \P X$ is a projective q-cofibration of $\topdgr_0$. The case where $\de Z \to Z$ is a \textit{trivial} projective q-cofibration is similar. 
\epf

We need for the proof of Theorem~\ref{PathSpacePreservesCofibrancy} some elementary information about the m-model structure and the h-model structure of $\top$ to complete the transfinite induction. We invite the reader to look \cite[Section~2]{leftproperflow} up.

\bth \label{PathSpacePreservesCofibrancy}
Let $X$ be a $\mathcal{P}$-flow. If $X$ is q-cofibrant, then the path $\mathcal{P}$-space functor $\P X$ is projective q-cofibrant. In particular, for every $(\alpha,\beta)\in X^0\p X^0$, the $\mathcal{P}$-space $\P_{\alpha,\beta}X$ is projective q-cofibrant if $X$ is q-cofibrant.
\eth

\bpf
By hypothesis, $X$ is q-cofibrant. Consequently, the map of $\mathcal{P}$-flows $\varnothing^\flat\to X$ is a retract of a transfinite composition $\liminj A_\lambda$ of pushouts of maps of the form $\globP(\de Z)\to \globP(Z)$ such that $\de Z \to Z$ is a projective q-cofibration of $\topdgr_0$ and of the maps $C:\varnothing^\flat \to \{0\}^\flat$ and $R:\{0,1\}^\flat\to \{0\}^\flat$. Since the set map $\varnothing^0=\varnothing \to X^0$ is one-to-one, one can suppose that the map $R:\{0,1\}^\flat\to \{0\}^\flat$ does not appear in the cellular decomposition. Moreover we can suppose that $A_0=X^0$ and that for all $\lambda\geq 0$, $A_\lambda \to A_{\lambda+1}$ is a pushout of a map of the form $\glob(\de Z)\to \glob(Z)$ where the map $\de Z\to Z$ is a projective q-cofibration of $\mathcal{P}$-spaces. 

Consider the set of ordinals $\{\lambda \mid \P A_\lambda \hbox{ not q-cofibrant}\}$. Note that this set does not contain $0$. If this set is nonempty, then it contains a smallest element $\mu_0>0$. By Theorem~\ref{pretower}, the ordinal $\mu_0$ is a limit ordinal. By definition of the ordinal $\mu_0$ and by Theorem~\ref{pretower} again, for all $\mu<\nu \leq \mu_0$, the maps of $\mathcal{P}$-spaces $\P A_{\mu} \to \P A_{\nu}$ are projective q-cofibrations of $\mathcal{P}$-spaces. By \cite[Proposition~7.1]{dgrtop}, we deduce that for all $\mu<\nu \leq \mu_0$, the maps of $\mathcal{P}$-spaces $\P A_{\mu} \to \P A_{\nu}$ are objectwise m-cofibrations of topological spaces, and therefore objectwise h-cofibrations of topological spaces. By \cite[Proposition~2.6]{leftproperflow}, we obtain that all $\mu<\nu \leq \mu_0$, the maps of $\mathcal{P}$-spaces $\P A_{\mu} \to \P A_{\nu}$ are objectwise relative-$T_1$ inclusions. By Theorem~\ref{global-path-almost-accessible}, we obtain that the canonical map 
\[\liminj_{\mu<\mu_0}\P A_\mu \to \P \liminj_{\mu<\mu_0} A_\mu = \P A_{\mu_0}\] is an objectwise homeomorphism. We deduce that $\P A_{\mu_0}$ is projective q-cofibrant: contradiction. 
\epf

\section{Comparison of flows and \texorpdfstring{$\mathcal{P}$}{Lg}-flows}
\label{equivalencePflow-flow}

\bd \cite[Definition~4.11]{model3} \label{def_flow}
A {\rm flow} is a $\mathbf{1}$-flow where $\mathbf{1}$ is the terminal category (see Example~\ref{terminal-reparam})). The corresponding category is denoted by $\dtop$. 
\ed

In other terms, a flow is a small semicategory enriched over the closed (semi)monoidal category $(\top,\p)$. Let us expand the definition above. A \textit{flow} $X$ consists of a topological space $\P X$ of execution paths, a discrete space $X^0$ of states, two continuous maps $s$ and $t$ from $\P X$ to $X^0$ called the source and target map respectively, and a continuous and associative map \[*:\{(x,y)\in \P X\p \P X; t(x)=s(y)\}\longrightarrow \P X\] such that $s(x*y)=s(x)$ and $t(x*y)=t(y)$.  A morphism of flows $f:X\longrightarrow Y$ consists of a set map $f^0:X^0\longrightarrow Y^0$ together with a continuous map $\P f:\P X\longrightarrow \P Y$ such that 
\[\begin{cases}
f^0(s(x))=s(\P f(x))\\
f^0(t(x))=t(\P f(x))\\
\P f(x*y)=\P f(x)*\P f(y).
\end{cases}\] Let \[\P_{\alpha,\beta}X = \{x\in \P X\mid s(x)=\alpha \hbox{ and } t(x)=\beta\}.\]

\begin{nota}
	The map $\P f:\P X\longrightarrow \P Y$ can be denoted by $f:\P X\to \P Y$ if there if no ambiguity. The set map $f^0:X^0\longrightarrow Y^0$ can be denoted by $f:X^0\longrightarrow Y^0$ if there if no ambiguity.
\end{nota}

The category $\dtop$ is locally presentable. Every set can be viewed as a flow with an empty path space. The obvious functor $\set \subset \dtop$ is limit-preserving and colimit-preserving. The following example of flow is important for the sequel: 

\begin{exa}
	For a topological space $Z$, let $\glob(Z)$ be the flow defined by 
	\[\begin{aligned}
	&\glob(Z)^0=\{0,1\}, \\ 
	&\P \glob(Z)= \P_{0,1} \glob(Z)=Z,\\ 
	&s=0,\  t=1.
	\end{aligned}\]
	This flow has no composition law.
\end{exa}

\begin{nota}
	\[\begin{aligned}
	&I^{gl} = \{\glob(\mathbf{S}^{n-1})\subset \glob(\mathbf{D}^{n}) \mid n\geq 0\}\\
	&J^{gl} = \{\glob(\mathbf{D}^{n}\p\{0\})\subset \glob(\mathbf{D}^{n+1}) \mid n\geq 0\}\\
	&C:\varnothing \to \{0\}\\
	&R:\{0,1\} \to \{0\}\\
	&\vI = \glob(\{0\})
	\end{aligned}\]
	where the maps of $I^{gl}$ are induced by the inclusions $\mathbf{S}^{n-1}\subset \mathbf{D}^{n}$ and the maps of $J^{gl}$ are induced by the mapping $(x_1,\dots,x_n)\mapsto (x_1,\dots,x_n,0)$. 
\end{nota}

The \textit{q-model structure} of flows $(\dtop)_q$ is the unique combinatorial model structure such that $I^{gl}\cup\{C,R\}$ is the set of generating q-cofibrations and such that $J^{gl}$ is the set of generating trivial q-cofibrations (e.g. \cite[Theorem~7.6]{QHMmodel} or \cite[Theorem~3.11]{leftdetflow}). The weak equivalences are the maps of flows $f:X\to Y$  inducing a bijection $f^0:X^0\iso Y^0$ and a weak homotopy equivalence $\P f:\P X \to \P Y$. The q-fibrations are the maps of flows $f:X\to Y$  inducing a q-fibration $\P f:\P X \to \P Y$ of topological spaces.

Let $X$ be a flow. The $\mathcal{P}$-flow $\moore(X)$ is the enriched semicategory defined as follows: 
\begin{itemize}[leftmargin=*]
	\item The set of states is $X^0$.
	\item The $\mathcal{P}$-space $\P_{\alpha,\beta}\moore(X)$ is the $\mathcal{P}$-space $\Delta_{\mathcal{P}^{op}}(\P_{\alpha,\beta}X)$.
	\item The composition law is defined, using Proposition~\ref{Ptenseur} as the composite map 
	\[
	\xymatrix@C=4em
	{
	\Delta_{\mathcal{P}^{op}}(\P_{\alpha,\beta}X)\ot  \Delta_{\mathcal{P}^{op}}(\P_{\beta,\gamma}X) \iso \Delta_{\mathcal{P}^{op}}(\P_{\alpha,\beta}X \p \P_{\beta,\gamma}X)   \fr{\Delta_{\mathcal{P}^{op}}(*)}& \Delta_{\mathcal{P}^{op}}(\P_{\alpha,\gamma} X).}
	\]
\end{itemize}

We obtain the 

\bp \label{def-preN}
The construction above yields a well-defined functor \[\moore:\dtop\to \dtopP.\]
\ep

Consider a $\mathcal{P}$-flow $Y$. For all $\alpha,\beta\in Y^0$, let $Y_{\alpha,\beta}=\liminj \P_{\alpha,\beta}Y$. Let $(\alpha,\beta,\gamma)$ be a triple of states of $Y$. By Proposition~\ref{tensor-product}, the composition law of the $\mathcal{P}$-flow $Y$ induces a continuous map \[Y_{\alpha,\beta} \p Y_{\beta,\gamma} \iso \liminj(\P_{\alpha,\beta}Y \ot \P_{\beta,\gamma}Y) \longrightarrow \liminj \P_{\alpha,\gamma}Y \iso Y_{\alpha,\gamma}\] which is associative. We obtain the 

\bp \label{def-N}
For any $\mathcal{P}$-flow $Y$, the data 
\begin{itemize}[leftmargin=*]
	\item The set of states is $Y^0$
	\item For all $\alpha,\beta\in Y^0$, let $Y_{\alpha,\beta}=\liminj \P_{\alpha,\beta}Y$
	\item For all $\alpha,\beta,\gamma\in Y^0$, the composition law $Y_{\alpha,\beta}\p Y_{\beta,\gamma}\to Y_{\alpha,\gamma}$
\end{itemize}
assemble to a flow denoted by $\lmoore(Y)$. It yields a well-defined functor \[\lmoore:\dtopP \to \dtop.\]
\ep

\bth \label{adj-cat-flow-Mooreflow}
There is an adjunction $\lmoore \dashv\moore$.
\eth

\bpf 
We have to prove that there is a natural bijection \[\dtopP(Y,\moore(X))\iso \dtop(\lmoore(Y),X)\] for all $\mathcal{P}$-flows $Y$ and all flows $X$. Since the functors $\lmoore:\dtopP\to \dtop$ and $\moore:\dtop\to\dtopP$ preserve the set of states, the only problem is to verify that everything is well-behaved with the path spaces, and more specifically with the composition laws. 

Let $f:\lmoore(Y)\to X$ be a map of flows. It induces for each pair of states $(\alpha,\beta)$ of $Y^0$ a map of topological spaces  $\P_{\alpha,\beta}\lmoore(Y)\to \P_{f(\alpha),f(\beta)}X$ such that for all triples $(\alpha,\beta,\gamma)$ of $Y^0$, the following diagram is commutative (the horizontal maps being the composition laws): 
\[
\xymatrix
{
	\P_{\alpha,\beta}\lmoore(Y) \p \P_{\beta,\gamma}\lmoore(Y) \fr{*}\fd{} & \P_{\alpha,\gamma}\lmoore(Y) \fd{}\\
	\P_{f(\alpha),f(\beta)}X \p \P_{f(\beta),f(\gamma)}X \fr{*} & \P_{f(\alpha),f(\gamma)}X.
}
\]
By adjunction, this means that the map of $\mathcal{P}$-spaces \[\P_{\alpha,\beta}X \longrightarrow \Delta_{\mathcal{P}^{op}}(\P_{f(\alpha),f(\beta)}X)\] is compatible with the composition laws, and then that we have a well-defined map of $\mathcal{P}$-flows from $Y$ to $\moore(X)$. 

Let $f:Y\to \moore(X)$ be a map of $\mathcal{P}$-flows. It induces for each pair of states $(\alpha,\beta)$ of $Y^0$ a map of $\mathcal{P}$-spaces $\P_{\alpha,\beta}Y\to \Delta_{\mathcal{P}^{op}}(\P_{f(\alpha),f(\beta)}X)$. We obtain by adjunction a map of topological spaces  $\P_{\alpha,\beta}\lmoore(Y)\to \P_{f(\alpha),f(\beta)}X$ and by naturality  of the adjunction the commutative diagram  (the horizontal maps being the composition laws): 
\[
\xymatrix
{
	\P_{\alpha,\beta}\lmoore(Y) \p \P_{\beta,\gamma}\lmoore(Y) \fr{*}\fd{} & \P_{\alpha,\gamma}\lmoore(Y) \fd{}\\
	\P_{f(\alpha),f(\beta)}X \p \P_{f(\beta),f(\gamma)}X \fr{*} & \P_{f(\alpha),f(\gamma)}X.
}
\]
\epf

Before ending the paper, we need to recall the

\bth (\cite[Theorem~7.6]{dgrtop}) \label{eq-topdgr-top} The Quillen adjunction 
\[\liminj \dashv \Delta_{\mathcal{P}^{op}} : [\mathcal{P}^{op},\top_q]^{proj}_0 \leftrightarrows \top_q\]
is a Quillen equivalence.
\eth

We can now conclude the first part of the proof of the existence of a zig-zag of Quillen equivalences between the q-model structure of multipointed $d$-spaces and the q-model structure of flows: 

\bth \label{eq1} The adjunction $\lmoore\dashv\moore$ of Theorem~\ref{adj-cat-flow-Mooreflow} is a Quillen equivalence between the q-model structure of $\dtopP$ and the q-model structure of $\dtop$.
\eth

\bpf We already know that $\moore:\dtop\to \dtopP$ is a right adjoint by Theorem~\ref{adj-cat-flow-Mooreflow}. It takes (trivial) q-fibrations of $\dtop$ to (trivial) q-fibrations of $\dtopP$ by Theorem~\ref{Moore-fibration} and Theorem~\ref{Moore-trivial-fibration}. Therefore it is a right Quillen adjoint.

Let $X$ be a flow (it is necessarily fibrant). The map of flows \[\lmoore(\moore(X)^{cof}) \longrightarrow X\] is bijective on states since the functors $\moore$ and $\lmoore$ preserves the set of states. For every pair $(\alpha,\beta)\in X^0\p X^0$, it induces the continuous map \[\liminj \P_{\alpha,\beta}\big(\moore(X)^{cof}\big) \longrightarrow \P_{\alpha,\beta}X.\]
By adjunction, we obtain the map of $\mathcal{P}$-spaces 
\[\P_{\alpha,\beta}\big(\moore(X)^{cof}\big) \longrightarrow \Delta_{\mathcal{P}^{op}}(\P_{\alpha,\beta}X) = \P_{\alpha,\beta}(\moore(X)).\]
The latter is a weak equivalence of the projective q-model structure of $\topdgr_0$ since the map $\moore(X)^{cof}\to \moore(X)$ is a weak equivalence of the q-model structure of $\dtopP$. By Theorem~\ref{PathSpacePreservesCofibrancy}, the $\mathcal{P}$-space $\P_{\alpha,\beta}\big(\moore(X)^{cof}\big)$ is projective q-cofibrant. It means that $\P_{\alpha,\beta}\big(\moore(X)^{cof}\big)$ is a cofibrant replacement of $\Delta_{\mathcal{P}^{op}}(\P_{\alpha,\beta}X)$ for the projective q-model structure of $\topdgr_0$. By Theorem~\ref{eq-topdgr-top}, we deduce that the continuous map \[\liminj \P_{\alpha,\beta}\big(\moore(X)^{cof}\big) \longrightarrow \P_{\alpha,\beta}X.\]
is a weak homotopy equivalence. In other terms, the left Quillen adjoint $\lmoore:\dtopP\to\dtop$ is homotopically surjective.

Let $Y$ be a q-cofibrant $\mathcal{P}$-flow. The map of $\mathcal{P}$-flows \[Y\longrightarrow \moore(\lmoore Y)\] is bijective on states since the functors $\moore$ and $\lmoore$ preserves the set of states. For every pair $(\alpha,\beta)\in X^0\p X^0$, the map $Y\to \moore(\lmoore Y)$ induces the map of $\mathcal{P}$-spaces 
\[
\P_{\alpha,\beta} Y \longrightarrow \Delta_{\mathcal{P}^{op}} (\liminj \P_{\alpha,\beta} Y).
\]
The $\mathcal{P}$-space $\P_{\alpha,\beta} Y$ is projective q-cofibrant by Theorem~\ref{PathSpacePreservesCofibrancy}. The above map is therefore an objectwise weak homotopy equivalence by Theorem~\ref{eq-topdgr-top}. In other terms, the left Quillen adjoint $\lmoore:\dtopP\to\dtop$ is a left Quillen equivalence. 
\epf

\appendix

\section{Changes made to the original paper \cite{Moore1}}
Here we provide details about the changes made to the original paper \cite{Moore1} in this erratum. We note that  the numbering of the theorems of the paper are unchanged.

\begin{enumerate}

\item  The abstract is slightly changed: the word ``symmetric'' and ``commutative'' are removed.

\item  Definition~\ref{def0}, Definition~\ref{def1}, Definition~\ref{def-reparam} are adapted.

\item  The following paragraph is added at the end of Section~\ref{param}: 

In the cases of $(\mathcal{G},+)$ and $(\mathcal{M},+)$, the functors $(\ell,\ell')\mapsto \ell+\ell'$ and $(\ell,\ell') \mapsto \ell'+\ell$ coincide on objects, but not on morphisms. The terminal category is a symmetric reparametrization category. We do not know if there exist symmetric reparametrization categories not equivalent to the terminal category.

\item Proposition~\ref{dec}, Lemma~\ref{petitcalcul}, Proposition~\ref{asso}, Theorem~\ref{closedsemi} are adapted: only half of the proofs is written, the proof of the ``symmetric'' statement is left to the reader; Notation~\ref{nota0} and Proposition~\ref{Ptenseur} are adapted as well.

\item  In the proof of Proposition~\ref{Ftenseur} and \ref{tensor-product}, some parameters $\ell$ are relocated because a reparame\-trization category is not supposed to be symmetric anymore.

\item  ``Commutative semigroup''  is replaced by ``semigroup'' (for the set of objects of a reparametrization category) everywhere in the paper.

\item  ``(Closed) symmetric semimonoidal structure'' is replaced by ``(biclosed) semimonoidal structure'' everywhere in the paper.

\item  First paragraph of Section~\ref{basic-property-moore-flow}: the word  ``symmetric'' is removed.

\item In the remark after Definition~\ref{exp}, ``the axioms of enriched categories are satisfied''` is replaced by ``the axioms of tensored and cotensored categories are satisfied''.

\item  The statement of Proposition~\ref{calculpushout} is changed: it is adapted from the statement for the closed symmetric monoidal category $(\top,\p)$; we have to modify the definition of $C_p(S)$ and of the map $C_p(i^T_S)$ to not change the order of the factors.

\end{enumerate}

\end{document}